\documentclass[12pt,a4paper]{article}
\usepackage[T1]{fontenc}
\usepackage{amsthm,amsmath,amssymb,graphics}
\usepackage{graphicx}
\usepackage{oldgerm}
\usepackage{psfrag}
\usepackage{calrsfs}
\newcommand{\field}[1]{\mathbb{#1}}\newcommand{\R}{\field{R}}

\newcommand{\Hp}{\field{H}}
\newcommand{\N}{\field{N}}

\newcommand{\Z}{\field{Z}}
\newcommand{\ch}{\mbox{ch}}

\newcommand{\eps}{\varepsilon}
\newcommand{\CR}{\mathcal{R}}
\newcommand{\CB}{\mathcal{B}}
\newcommand{\ML}{\mathcal{ML}}
\theoremstyle{plain}
\newtheorem*{theorem}{Theorem 1}

\newtheorem*{theoremc}{Theorem 1 corrected}
\newtheorem{prop}{Proposition}[section]
\newtheorem{lemme}[prop]{Lemma}

\newtheorem{claim}[prop]{Claim}

\begin{document}

\title{Continuity of the bending map}
\author{Cyril Lecuire}
\date{}

\maketitle
{\bf R\'esum\'e}\\
L'application de plissage d'une vari\'et\'e hyperbolique de dimension $3$ associe \`a une m\'etrique hyperbolique convexe cocompacte sur une vari\'et\'e compacte \`a bord sa lamination g\'eod\'esique mesur\'ee de plissage. Il a \'et\'e d\'emontr\'e dans \cite{kes} et \cite{kami} que cette application est continue. Dans ce texte, on \'etudie l'extension de cette application \`a l'espace des m\'etriques hyperboliques g\'eom\'etriquement finies. On introduit une relation d'\'equivalence dans l'espace des laminations g\'eod\'esiques mesur\'ees et on montre que l'application quotient de l'application de plissage est continue.\\

{\bf Abstract}\\
The bending map of a hyperbolic $3$-manifold  maps a convex cocompact hyperbolic metric on a $3$-manifold with boundary to its bending measured geodesic lamination. As proved in \cite{kes} and \cite{kami}, this map is continuous. In the present paper we study the extension of this map to the space of geometrically finite hyperbolic metrics. We introduce a relationship on the space of measured geodesic laminations and shows that the quotient map obtained from the bending map is continuous. 

\newpage

\section*{Introduction}

\indent
Let $M$ be a compact, orientable $3$-manifold with boundary.  Assume that $M$ is hyperbolic, namely that the interior of $M$ is endowed with a complete metric $\sigma$ of constant sectional curvature $-1$. Assume also that $\partial M$ contains a surface with genus greater than $1$. A fundamental subset of $(M,\sigma)$ is its convex core $N(\sigma)$. This core $N(\sigma)$ is defined as the smallest non-empty closed subset of the interior of $M$ which is locally convex and homotopically equivalent to $M$. Its boundary $\partial N(\sigma)$ endowed with the intrinsic metric (given by rectifiable path length) is isometric to a hyperbolic surface of finite volume and can be embedded in a natural way into $\partial M$. This surface is bent along a geodesic lamination and the amount of bending is described by a measured geodesic lamination called the {\em bending measured geodesic lamination} of $\sigma$ (cf. \cite{notes} or \cite{ceg}). This yields a bending map $b$ which to a complete hyperbolic metric associates its bending measured geodesic lamination.\\
\indent
In \cite{bodiff}, F. Bonahon considers quasi-isometric deformations of a given metric $\sigma$ on $int(M)$, namely hyperbolic metrics $\sigma'$ on $int(M)$ for which there exists a diffeomorphism $(int(M),\sigma)\rightarrow(int(M),\sigma')$ whose differential is uniformly bounded. Let ${\cal QD}(\sigma)$ be the space of quasi-isometric deformations of a given metric $\sigma$, where we identify two deformations $(int(M),\sigma)\rightarrow(int(M),\sigma')$ and $(int(M),\sigma)\rightarrow(int(M),\sigma'')$ if they are isotopic. The continuity of the bending map $b:{\cal QD}(\sigma)\rightarrow {\cal ML}(\partial M)$ is proved in \cite{kami}, using ideas of Thurston. Its differentiability (in a weak sense) is proved in \cite{bodiff}.\\
\indent
A complete hyperbolic metric $\sigma$ on $int(M)$ is said to be {\em convex cocompact} if $N(\sigma)$ is compact. If $\sigma$ is convex cocompact, then ${\cal QD}(\sigma)$ is the set of convex cocompact metrics on $int(M)$. In this case, the continuity of $b:{\cal QD}(\sigma)\rightarrow {\cal ML}(\partial M)$ has been proved in \cite{kes} and its image has been described in \cite{meister1} and \cite{espoir}. A complete hyperbolic metric $\sigma$ on $int(M)$ is said to be {\em geometrically finite} if $N(\sigma)$ has finite volume. When $\sigma$ is a geometrically finite metric, ${\cal QD}(\sigma)$ is the set of geometrically finite metrics having the same parabolic subgroups as $\sigma$. In the present paper, we are addressing the question of the continuity of the bending map on the whole set of geometrically finite metrics. In particular, we are interested in sequence converging to a limit with some new parabolics.\\
\indent
Since we want to consider metrics which do not have the same parabolic subgroups, we have to allow deformations which are not quasi-isometric. We will consider the set of isotopy classes of hyperbolic metrics on the interior of $M$. Two metrics $\sigma_1$ and $\sigma_2$ are identified if there exists a diffeomorphism $f : M\rightarrow M$ isotopic to the identity such that $\sigma_1=f^*\sigma_2$. We will consider the set ${\cal GF}(M)$ of isotopy classes of geometrically finite hyperbolic metrics which are not fuchsian. We will topologise ${\cal GF}(M)$ in the following way. Let us choose a point $x$ in $int(M)$. A metric $\sigma_2$ lies in a $(k,r)$-neighbourhood of $\sigma_1$ if there exists a diffeomorphism $g:M\rightarrow M$ isotopic to the identity such that the restriction of $g$ to the ball $B(x,r)\subset (M,\sigma_1)$ is a $k$-quasi-isometry into its image in $(M,\sigma_2)$. We obtain a basis of neighbourhoods of $\sigma_1$ by letting $k$ tend to $1$ and $r$ tend to $+\infty$. The topology defined in this way does not depend on the choice of the point $x$. For a metric $\sigma\in{\cal GF}(M)$, ${\cal QD}(\sigma)$ can be viewed as a subset of ${\cal GF}(M)$. The topology of ${\cal QD}(\sigma)$ considered as a subset of ${\cal GF}(M)$ coincides with the topology given by quasi-isometric deformations\\

\indent
The bending map $b_{\cal GF}(M):{\cal GF}(M)\rightarrow {\cal ML}(\partial M)$ maps an isotopy class of geometrically finite metrics to its bending measured geodesic lamination. The image of $b_{\cal GF}$ has been described in \cite{meister1} and \cite{espoir}, it is the set ${\cal P}(M)$ of measured geodesic laminations satisfying the following conditions:
\begin{description}
\item -  a) no closed leaf of $\lambda$ has a weight greater than $\pi$;
\item -  $b)$ $\exists\eta>0$ such that, for any essential annulus $E$, $i(\partial E,\lambda)\geq\eta$;
\item -  c) $i(\lambda,\partial D)>2\pi$ for any essential disc $D$.
\end{description}

\indent
Taking a careful look at the behaviour of the map $b_{\cal GF}$, we notice that it is not a continuous map. If a metric $\sigma$ lies in ${\cal GF}(M)$ but has some rank one cusps, its bending measured geodesic lamination $\lambda$ has some compact leaves with a weight equal to $\pi$. Let us denote by $\lambda^{(p)}$ the union of the leaves of $\lambda$ which have a weight equal to $\pi$. Using the result of \cite{meister1}, it is not hard to construct a sequence of metrics $\sigma_n\in{\cal GF}(M)$ with measured geodesic laminations $\lambda_n$ such that the sequence $(\lambda_n)$ converges to a measured geodesic lamination $\lambda_{\infty}$ which differs from $\lambda$ only on $\lambda^{(p)}$ and which has some leaves with a weight greater than $\pi$. Using some arguments of \cite{espoir}, we get that a subsequence of $(\sigma_n)$ converges to a geometrically finite metric $\sigma_\infty$. Since $\lambda_\infty$ does not satisfies condition $a)$, it is not the bending measured geodesic lamination of $\sigma_\infty$. Thus we get that $b_{\cal GF}$ is not continuous on any neighbourhood of a metric with some rank one cusps. To overcome this difficulty, we will quotient the space ${\cal ML}(\partial M)$ of measured geodesic laminations by the following relationship :\medskip\\
\indent
Let  $\lambda, \mu \in {\cal ML}(\partial M)$ be two measured geodesic laminations and let us denote by $\lambda'$ (resp. $\mu'$) the measured geodesic laminations obtained  by replacing  by $\pi$ the weights of the leaves of $\lambda$ (resp. $\mu$) which have a weight greater than $\pi$; we set  $\lambda {\cal R} \mu$ if and only if $\lambda'=\mu'$. We denote by $\dot\lambda$ the class of $\lambda$ modulo ${\cal R}$.\medskip\\
\indent Let us endow ${\cal ML}(\partial M)$ with the weak$^*$ topology and ${\cal ML}(\partial M)/{\cal R }$ with the quotient topology. From $b_{\cal GF}(M)$ we obtain a quotient map $b_{\cal R} :{\cal GF}(M)\rightarrow {\cal P}(M)/{\cal R}$.  We will prove the following result :

\begin{theorem}         \label{two}
The map $b_{\cal R}$ from ${\cal GF}(M)$ to ${\cal ML}(\partial M)/{\cal R}$ is a continuous map.
\end{theorem}

\indent
In \cite{pasf}, we show the reverse of this theorem. This gives rise to a criterion for the strong convergence of a sequence of geometrically finite representations $\rho_n\in U(M)$.\\

\indent The paper is organised as follows. In section 1, we state some definitions and we prove some facts about ${\cal ML}(\partial M)/{\cal R}$. In section 2, we study convex pleated surfaces and prove the continuity of the bending measured geodesic laminations of a converging sequence of convex pleated surfaces. In section 3 we use the results of section 2 to prove the continuity of $b_{\cal R}$.

\section{Definitions}
\indent Let $\sigma$ be a hyperbolic metric (up to isotopy) on $int(M)$. Given an isometry from the interior of $\tilde M$ to $\Hp^3$, the covering transformations yield a discrete faithful representation\linebreak $\rho:\pi_1(M)\rightarrow Isom(\Hp^3)$. The representations that appear in this way will be called {\em representations associated to $\sigma$}. The set of representations associated to $\sigma$ is the set of all representations conjugated to $\rho$. The image  $\rho(\pi_1(M))$ is a finitely generated torsion free Kleinian group and $int(M)$ endowed with $\sigma$ is isometric to $\Hp^3/\rho(\pi_1(M))$. The {\em Nielsen core} of $\Hp^3/\rho(\pi_1(M))$ is the quotient by $\rho(\pi_1(M))$ of the convex hull $C(\rho)$ of the limit set $L_{\rho}$ of $\rho(\pi_1(M))$ (see \cite[chap 8]{notes} for details). This set $N(\rho)$ and the convex core defined in the introduction are isometric and from now on we will identify them. The {\em thick part}, $N(\rho)^{ep}$ of the Nielsen core is the complementary of the cuspidal part of $\Hp^3/\rho(\pi_1(M))$ in $N(\rho)$. The representation $\rho$ is {\em geometrically finite} when $N(\rho)^{ep}$ is compact (here it is equivalent to say that $N(\rho)$ has finite volume) and {\em convex cocompact} when $N(\rho)$ is compact.\\
\indent
When $\rho$ is geometrically finite and not fuchsian, the natural retraction from\linebreak $\Hp^3/\rho(\pi_1(M))$ to $N(\rho)$ associates to $\sigma$ a homeomorphism (defined up to isotopy)\linebreak $h:M\rightarrow N(\rho)^{ep}$. Such a homeomorphism will be said to be {\em associated to $\sigma$}.\\
\indent Let $\sigma\in{\cal QD}(M)$ be a metric whose only cusps are rank $2$ cusps and let $\rho$ and $h$ be respectively a representation and a homeomorphism associated to $\sigma$. Let $\tilde h:\tilde M\rightarrow\Hp^3$ be a lift of $h$, we will define the Floyd-Gromov compactification $\overline{\tilde M}$ of $\tilde M$ as being the closure of $\tilde h(\tilde M)$ in the usual compactification of $\Hp^3$ by the unit ball.\\
\indent Let $(\sigma_n)$ be a sequence of isotopy classes of complete hyperbolic metrics on the interior of $M$. The sequence $(\sigma_n)$ {\em converges algebraically} when there is a sequence of representations $\rho_n:\pi_1(M)\rightarrow Isom(\Hp^3)$ associated to $\sigma_n$ (as above) that converges algebraically; namely $\rho_n(g)$ converges for any $g\in\pi_1(M)$. We obtain a new representation $\rho_{\infty}:\pi_1(M)\rightarrow Isom(\Hp^3)$ defined by $\rho_{\infty}(g)=\lim_{n\longrightarrow\infty}\rho_n(g)$ for any $g\in\pi_1(M)$. This representation $\rho_\infty$ is discrete and faithful (cf. \cite{jorgi}). This representation defines a metric on a manifold homotopy equivalent to $M$. This manifold might not be homeomorphic to $M$ (examples are given in \cite{cancan}).\\
\indent
The sequence $(\sigma_n)$ {\em converges geometrically} if there is a sequence of representations\linebreak $\rho_n:\pi_1(M)\rightarrow Isom(\Hp^3)$ associated to the $\sigma_n$ such that $(\rho_n(\pi_1(M))$ converges geometrically. The sequence of groups  $(\rho_n(\pi_1(M))$ converges geometrically to a group  $\Gamma_\infty\subset Isom(\Hp^3)$ if and only if :\\
\indent - for any sequence $a_n\in\rho_n(\pi_1(M))$, any accumulation point $a_\infty$ of $\{a_n|\,n\in\N\}$ lies in $\Gamma_\infty$;\\
\indent - any element $a_\infty$ of $\Gamma_\infty$ is the limit of a sequence $a_n\in\Gamma_n$.\\
\indent 
The sequence $\rho_n(\pi_1(M))$ {\em converges strongly} if there is a sequence of representations $\rho_n$ associated to $\sigma_n$ such that $(\rho_n)$ converges algebraically to a representation $\rho_\infty$ and that $(\rho_n(\pi_1(M)))$ converges geometrically to $\rho_\infty(\pi_1(M))$. If $(\sigma_n)$ converges to $\sigma$ for the topology defined in the introduction, then  $(\sigma_n)$ converges strongly to $\sigma_\infty$ (cf. \cite{ceg}).

\section{Geodesic laminations and the relationship ${\cal R}$}
\indent A  {\em geodesic lamination} $L$ on $\Hp^2$ is a closed subset which is the disjoint union of complete geodesics. A complete geodesic lying in $L$ is a {\it leaf} of $L$.\\
\indent
A  {\em measured geodesic lamination} $\lambda$ is a transverse measure for some geodesic lamination $|\lambda|$. Any arc $k\approx [0,1]$ embedded in $S$ transversely to $|\lambda|$, such that $\partial k\subset S-\lambda$, is endowed with an  additive measure $d\lambda$ such that :\\
\indent - the support of $d\lambda_{|k}$ is $|\lambda|\cap k$;\\
\indent - if an arc $k$ can be homotoped into $k'$ by a homotopy respecting $|\lambda|$ then $\int_k\! d\lambda=\int_{k'} d\lambda$.\\
We will denote by ${\cal ML}(\Hp^2)$ the space of measured geodesic lamination topologised with the topology of the weak$^*$ convergence.\\

\indent 
Let $S$ be a surface (which may not be compact) endowed with a complete hyperbolic metric with finite area. A geodesic lamination in $S$ is a compact subset which is the disjoint union of simple complete geodesics.  Using the fact that two complete hyperbolic metrics with finite area on $S$ are quasi-isometric, this definition can be made independent of the chosen metric on $S$ (see for example \cite{conti}). A measured geodesic lamination is a transverse measure for some geodesic lamination as defined above. Let $\gamma$ be a weighted simple closed geodesic with support $|\gamma|$ and weight $w$ and let $\lambda$ be a measured geodesic lamination. The intersection number of $\gamma$ and $\lambda$ is defined by $i(\gamma,\lambda)=w \int_{|\gamma|} d\lambda$. If $|\gamma|$ is a leaf of $\lambda$, then we define the intersection number by $i(\gamma,\lambda)=0$. The weighted simple closed curves are dense in ${\cal ML}(S)$ and this intersection number extends continuously to a function $i:{\cal ML}(S)\times{\cal ML}(S)\rightarrow\R$ (cf. \cite{bouts}).\\
\indent
A measured geodesic lamination $\lambda\in{\cal ML}(S)$ is {\em arational} if for any essential simple closed curve $c$ which is not homotopic to a cusp, we have $i(c,\lambda)=\int_c d\lambda>0$.\\

\indent Let us recall the definition of  ${\cal R}$ given in the introduction.\smallskip

\indent Let  $\lambda, \mu \in {\cal ML}( S)$ be two measured geodesic laminations. Let us denote by $\lambda'$ (resp. $\mu'$) the measured geodesic lamination obtained by replacing by $\pi$ the weights of the compact leaves of $\lambda$ (resp. $\mu$) which have a weight greater than $\pi$. We will say that $\lambda$ is related to $\mu$ by the relationship ${\cal R}$, and we will write $\lambda {\cal R} \mu$, if and only if $\lambda'=\mu'$. We will denote by $\dot\lambda$ the class of $\lambda$ modulo ${\cal R}$ and we will topologise ${\cal ML}(S)$ with the quotient topology of the weak$^*$ topology on ${\cal ML}(S)$.\\
\indent
Let $\partial_{\chi<0} M$ be the union of the connected components of $\partial M$ with negative Euler characteristic. To simplify the notations, we will note ${\cal ML}(\partial M)$ for ${\cal ML}(\partial_{\chi<0} M)$.\medskip\\

\indent
Let $\lambda,\mu\in{\cal ML}(\partial M)$ be two measured geodesic laminations. If $\lambda{\cal R}\mu$, then $\lambda$ and $\mu$ share the same support. Thus we can define the support $|\dot\lambda|$ of an element $\dot\lambda$ of ${\cal ML}(\partial M)$ as being the support of any representative of $\dot\lambda$.\\
\indent Let $(\lambda_n)$ be a sequence of measured geodesic laminations such that $(\dot\lambda_n)$ converges to $\dot\lambda$ in ${\cal ML}(\partial M)/{\cal R}$ and let $\lambda\in{\cal ML}(\partial M)$ be a representative of $\dot\lambda$. Let us denote by $\lambda^{(p)}$ the union of the compact leaves of $\lambda$ with a weight at least $\pi$. We have the following :

\begin{claim}   \label{nini}
Let $\lambda\in{\cal ML}(\partial M)$ and let $\lambda'$ be the representative of $\dot\lambda$ whose compact leaves have all a weight at most $\pi$. If $k\subset \partial M$ is a simple arc such that the points of $k\cap|\dot\lambda|$ and of $k\cap|\dot\lambda_n|$ are transverse intersections, then we have $\int_k\! d\lambda'\leq\liminf\int_k\! d\lambda_n$. Furthermore, if $k$ does not intersect $\lambda^{(p)}$, then $\int_k\! d\lambda_n$  converges to $\int_k\! d\lambda$. 
\end{claim}

\begin{proof}
Let $k\subset\partial M$ be a simple arc such that the points of $k\cap|\dot\lambda|$ and of $k\cap|\dot\lambda_n|$ are transverse intersections. By definition of $\lambda'$,  up to cutting $k$ into finitely many sub-arcs, we may assume that we have $\int_k d\lambda'\leq\pi$.  The set ${\cal V}_{k,\eps}(\lambda)=\{\gamma\in{\cal ML}(\partial M)/ |\int_k\! d\gamma-\int_k\! d\lambda|<\eps \}$ is a neighbourhood of $\lambda$ in ${\cal ML}(\partial M)$. Since $\dot\lambda_n$ converges to $\dot\lambda$, for any $\eps$, there is $n_\eps$ such that for any $n\geq n_\eps$, $\exists \mu_n,\alpha_n$ with $\mu_n{\cal R}\lambda$, $\alpha_n{\cal R}\lambda_n$ and $\alpha_n\in {\cal V}_{k,\eps}(\mu_n)$. So, for $n\geq n_\eps$, we have $|\int_k\! d\alpha_n-\int_k\! d\mu_n|<\eps$. Since $\lambda'$ is the representative of $\dot\lambda$ whose compact leaves have all a weight at most $\pi$ the measure of any $\mu\in{\cal ML(\partial M)}$ satisfying $\mu{\cal R}\lambda$ is at least the measure of $\lambda'$, namely we have $\int_k\! d\lambda'\leq\int_k\! d\mu$. It follows that we have $\int_k\! d\alpha_n\geq\int_k\! d\mu_n-\eps\geq\int_k\! d\lambda'-\eps$. Thus we get $\int_k\! d\lambda'\leq\liminf\int_k\! d\alpha_n$.\\
\indent
If we have $\int_k\! d\lambda_n<\pi$, then $k$ does not intersect any  closed leaf of $\lambda_n$ with a weight at least $\pi$. Therefore we have $\int_k\! d\lambda_n=\int_k\! d\alpha_n$ and the inequality $\int_k\! d\lambda'\leq\liminf\int_k\! d\lambda_n$ follows from the paragraph above. Otherwise we have $\int_k\! d\lambda_n>\pi\geq \int_k\! d\lambda'$ by assumption and the inequality is obvious.\\
\indent
If $k$  does not intersect any closed leaf of  $\lambda'$ with a weight equal to $\pi$, up to cutting $k$ into finitely many sub-arcs, we may assume that we have $\int_k d\lambda'<\pi$. Since $k$ does not intersect any closed leaf of  $\lambda'$ with a weight equal to $\pi$, we have $\int_k\! d\mu=\int_k\! d\lambda$ for any $\mu\in{\cal ML}(\partial M)$ satisfying $\mu{\cal R}\lambda$.  Especially we have $\int_k\! d\alpha_n\longrightarrow\int_k\! d\lambda'<\pi$. Hence we have $\int_k\! d\alpha_n<\pi$ for $n$ large enough. It follows that $k$ does not intersect any closed leaf of $\alpha_n$ with a weight at least $\pi$. So we have $\int_k\! d\lambda_n=\int_k\! d\alpha_n\longrightarrow\int_k\! d\lambda'=\int_k\! d\lambda$.
\end{proof}

\indent
An arc $k$ is {\it generic} if it is transverse to every simple geodesic of $S$. Especially a generic arc is transverse to every geodesic lamination. By \cite{birse}, the union of all simple geodesics of $S$ has Hausdorff dimension $1$. It follows that almost every geodesic arc is generic and that every arc can be approximated by a generic arc.

\begin{claim} \label{cpir}
Let $(\lambda_n)\in{\cal ML}(\partial M)$ be a sequence of measured geodesic laminations such that $(\dot\lambda_n)$ converges to $\dot\lambda$ and that $|\lambda_n|$ converges to some geodesic lamination $L $ in the Hausdorff topology. We have $|\dot\lambda|\subset L $.
\end{claim}

\begin{proof}
Let $x$ be a point of $|\dot\lambda|$, let $\eps>0$ be a real number and let $k$ be a geodesic generic arc intersecting $|\dot\lambda|$, with length $\eps$ such that $x$ lies in the interior of $k$. Since we have $\int_{k} d\lambda'>0$, we deduce from \ref{nini} that, for $n$ large enough, we have $\int_k\! d\lambda_n\geq\frac{\int_k\! d\lambda'}{2}>0$. Therefore $k$ intersects $\lambda_n$ and $x$ lies in an $\eps$-neighbourhood of $|\lambda_n|$. Considering a covering of $|\dot\lambda|$ by discs with radius $\eps$  and with centres lying in $|\dot\lambda|$, we get that, for $n$ large enough, $|\dot\lambda|$ lies in an $\eps$-neighbourhood of $|\lambda_n|$. Letting $\eps$ tend to $0$, we get $|\dot\lambda|\subset L $.
\end{proof}

\indent
Claim  \ref{nini} can also be used to prove that the space ${\cal ML}(\partial M)/{\cal R}$ is a Hausdorff space.

\begin{lemme}
The space ${\cal ML}(\partial M)/{\cal R}$ is a Hausdorff space.
\end{lemme}

\begin{proof}
Let $\dot\lambda$ and $\dot\mu$ be two elements of ${\cal ML}(\partial M)/{\cal R}$ such that any neighbourhood of $\dot\lambda$ intersects any neighbourhood of $\dot\mu$. So there is a sequence of measured geodesic laminations $\lambda_n$ such that $\dot\lambda_n$ converges simultaneously to $\dot\lambda$ and to $\dot\mu$. Let us denote by $\lambda^{(p)}$ (resp. $\mu^{(p)}$) the union of the compact leaves of $\lambda$ (resp. $\mu$) with a weight at least $\pi$. Let $\lambda'$ (resp. $\mu'$) be the representative of $\dot\lambda$ (resp. $\dot\mu$) whose compact leaves have all a weight at most $\pi$. Let $k\subset \partial M-\lambda^{(p)}$ be a generic arc intersecting $|\lambda|$ and $|\mu|$ transversely so that $\int_k\! d\lambda<\pi$. By Claim \ref{nini}, we have $\int_k\! d\mu'= \lim \int_k\! d\lambda_n<\pi$. It follows that $\mu^{(p)}\subset\lambda^{(p)}$. Reversing the roles of  $\lambda$ and of $\mu$, we get $\mu^{(p)}=\lambda^{(p)}$. It follows also from Claim \ref{nini} that we have the equality $\int_k\! d\mu=\int_k\! d\lambda$ for any arc $k\subset \partial M-\lambda^{(p)}$. This yields the conclusion $\dot\lambda=\dot\mu$.
\end{proof}

\indent
The following variation of Claim \ref{nini} will be used in the present paper.

\begin{claim}   \label{ni}
If $c$ is a simple closed curve that does not intersect $\lambda^{(p)}$ transversely, then the sequence $i(c,\lambda_n)$ converges to $i(c,\lambda)$.
\end{claim}

\begin{proof}
If the points of $c$ intersects $|\dot\lambda|$ transversely, then by Claim \ref{cpir}, $c$ intersects $|\dot\lambda|$ transversely for $n$ large enough. In this case, we get the conclusion by cutting $c$ into two arcs and by applying Claim \ref{nini} to these two arcs.\\
\indent Let us now consider the case where $c$ does not intersect $|\dot\lambda|$ transversely. Let $\lambda\in{\cal ML}(S)$ be a representative of $\dot\lambda$. Consider a simple closed curve $c'$ which is disjoint from $|\lambda|$, which is freely homotopic to $c$ and which is the union of $2$ generic arcs. We have $\int_{c'}\! d\lambda=0$. By Claim \ref{nini}, the sequence $\int_{c'}\! d\lambda_n$ converges to $\int_{c'}\! d\lambda=0$. Furthermore for any measured geodesic lamination $\gamma\in{\cal ML}(S)$, we have $i(c,\gamma)\leq \int_{c'}\! d\gamma$. Thus we can conclude that $i(c,\lambda_n)$ converges to $0=i(c,\lambda)$.
\end{proof}

\section{Convex pleated surfaces}
A {\em pleated surface} in a complete hyperbolic $3$-manifold $M$ is a map $f:S\rightarrow M$  from a surface $S$ to $M$ with the following properties :
\begin{description}
\item - the path metric obtained by pulling back the hyperbolic metric of $M$ by $f$ is a hyperbolic metric $s$ on $S$;
\item - every point of $S$ lies in the interior of some $s$-geodesic arc which is mapped into a geodesic arc in $M$;
\item - if $c\subset S$ is a simple closed curve lying in a cusp of $S$ and if $c$ does not bound a disc in $S$, then $f(c)$ does not bound a disc in $M$.
\end{description}

\indent
A map $\hat f:\Hp^2\rightarrow\Hp^3$ is a pleated map if any point of $\Hp^2$ lies in the interior of a geodesic arc which is mapped by $\hat f$ into a geodesic arc.\\
\indent
The {\em pleating locus} of a pleated map is the set of points of $\Hp^2$ where the map fails to be an isometry. The pleating locus of a pleated map is a geodesic lamination (cf. \cite{notes}).\\

{\em An abstract pleated surface} is a triple $(\hat f,\Gamma,\rho)$ where $\hat f:\Hp^2\rightarrow\Hp^3$ is a pleated map, $\Gamma$ is a lattice in $Isom(\Hp^2)$ and $\rho:\Gamma\rightarrow Isom(\Hp^3)$ is a discrete representation (which may not be faithful) such that for any $a\in\Gamma$, we have $\hat f\circ a=\rho(a)\circ \hat f$ and that if $a\in\Gamma$ is a parabolic isometry then $\rho(a)$ is also a parabolic isometry. 

\indent Abstract pleated surfaces and pleated surfaces are related as follows:\\
\indent
When $\rho(\Gamma)$ has no torsion, the abstract pleated surface $(\hat f,\Gamma,\rho)$ induces a pleated surface $f:S\rightarrow M$ where $S\approx\Hp^2/\Gamma$, $M\approx\Hp^3/\rho(\Gamma)$ and where $f$ is the quotient map from $\hat f$.\\
\indent If $f:S\rightarrow M$ is a pleated surface, consider isometric covering maps $\Hp^2\rightarrow S$ and\linebreak $\Hp^3\rightarrow M$. These maps yield representations $r:\pi_1(S)\rightarrow Isom(\Hp^2)$ and\linebreak $R:\pi_1(M)\rightarrow Isom(\Hp^3)$ and by lifting $f$ to a map $\hat f:\Hp^2\rightarrow \Hp^3$ we get an abstract pleated surface $(\hat f,r(\pi_1(S)),R\circ f_*)$.\\
\indent
In the following we will omit the adjective abstract and assume that our (abstract) pleated surfaces are torsion free.\\

\indent
We will consider the following topology on the space of abstract pleated surfaces.\\
A sequence $(\hat f_n,\Gamma_n,\rho_n)$ of pleated surfaces converge to a pleated surface $(\hat f,\Gamma,\rho)$ if and only if :
\begin{description}
\item - $(\Gamma_n)$ converges geometrically to $\Gamma$;
\item - for any sequence $a_n\in\Gamma_n$ converging to $a\in\Gamma$, $\rho_n(a_n)$ converges to $\rho(a)$;
\item - $\hat f_n$ converges to $\hat f$ on any compact set of $\Hp^2$.
\end{description}

\indent
Let $(\hat f,\Gamma,\rho)$ be a (abstract) pleated surface (without torsion), let $\hat L \subset\Hp^2$ be the pleating locus of $\hat f$ and let $P$ be a connected component of $\Hp^2-\hat L$. The surface $\hat f(P)$ lies in a geodesic plane $\Pi_P$. Given an orientation of $\Hp^2$, $\Pi_P$ inherits a natural orientation and we denote by $H^+_P$ (resp. $H^-_P$) the half-space bounded by $\Pi_P$ such that the union of a direct frame of $\Pi_P$ and of the inward normal vector to $\Pi_P$ is a direct frame of $\Hp^3$ (resp. indirect).\\

A pleated surface $(\hat f,\Gamma,\rho)$ with pleating locus $\hat L $ is a {\em convex pleated surface} if :
\begin{description}
\item 1) there is $\epsilon\in\{+,-\}$ such that for any component $P$ of $\Hp^2-\hat L $, $\hat f(\Hp^2)$ lies in $H^\epsilon_P$;
\item 2) the interior of $C_{\hat f}=\bigcap\{H_{P_i}^{\epsilon}|P_i$ is a connected component of $\Hp^2-\hat L\}$ is not empty.
\end{description}

\indent If a pleated surface $(\hat f,\Gamma,\rho)$ satisfies $1)$ but not $2)$, namely if $C_{\hat f}$ as empty interior, we will call it an {\em even pleated surface}.

\indent
In the following, for any convex pleated surface $(\hat f,\Gamma,\rho)$, we will choose the orientation of $\Hp^2$ so that for any component $P$ of $\Hp^2- L$, $\hat f(\Hp^2)$ lies in $H^+_P$.

\begin{lemme}   \label{confer}
The set of the pleated surfaces which are either convex or even is a closed subset of the set of pleated surfaces.
\end{lemme}

\begin{proof}
Let $(\hat f_n,\Gamma_n,\rho_n)$ be a sequence of convex and even pleated surfaces converging to\linebreak $(\hat f_\infty,\Gamma_\infty,\rho_\infty)$. 
 Let us follow \cite[Lemmas 20 and 21]{meister1}, to show that $\hat f_\infty$ is either a convex pleated surface or an even pleated surface. Let us denote by $\hat L_n$ the pleating locus of $(\hat f_n,\Gamma_n,\rho_n)$ and let us consider a geodesic lamination $\hat L_\infty$ which is a limit point of $\{\hat L_n/n\in\N\}$ in the Hausdorff topology. By \cite{ceg} the pleating locus of $\hat f_\infty$ lies in $\hat L_\infty$. Let us consider a component $P_\infty$ of $\Hp^2-\hat L_\infty$ and a component $P_n$ of $\Hp^2-\hat L_n$ such that $(P_n)$ tends to $P_\infty$. Since $\hat f_n(P_n)$ converges to $\hat f_\infty(P_\infty)$, up to extracting a subsequence, $H^+_{P_n}$ converges to a half-space $H^+_{P_\infty}$ such that $\hat f_\infty(P_\infty)\subset \partial H^+_{P_\infty}$. Since $\hat f_n(\Hp^2)$ converges to $\hat f_\infty(\Hp^2)$, we have $\hat f_\infty(\Hp^2)\subset H^+_{P_\infty}$. Doing this for any component of $\Hp^2-L_\infty$, we conclude that $\hat f_\infty$ satisfies $1)$.
\end{proof}

\indent
Let $\hat x$ be a point of $\Hp^2$, a {\em support plane} of $\hat f(\Hp^2)$ at $\hat f(\hat x)$ is a hyperbolic plane $\Pi_{\hat f(\hat x)}$ containing $\hat f(\hat x)$ and such that $\hat f(\Hp^2)$ lies entirely in one of the two half-spaces bounded by $\Pi_{\hat f(\hat x)}$. We will denote this half-space by $H^+_{\hat f(\hat x)}$. Let $\hat k\subset\Hp^2$ be a compact geodesic segment, a {\em polygonal approximation to $\hat f(\hat k)$} is a finite family ${\textgoth P}=\{(\hat x_i,\Pi_{\hat f(\hat x_i)})| i=1,\ldots, p\}$ such that :
\begin{description}
\item - $\partial\hat k=\{\hat x_1,\hat x_p\}$;
\item - the $\hat x_i$ are ordered points of $\hat k$;
\item - $\Pi_{\hat f(\hat x_i)}$ is a support plane at $\hat f(\hat x_i)$;
\item - $\Pi_{\hat f(\hat x_i)}\cap\Pi_{\hat f(\hat x_{i+1})}\neq\emptyset$ for any $i=1,\ldots, p-1$;
\item - if $\hat x_{i-1}=\hat x_i=\hat x_{i+1}$, then either $\Pi_{\hat f(\hat x_i)}=\Pi_{\hat f(\hat x_{i-1})}$ or $\Pi_{\hat f(\hat x_i)}=\Pi_{\hat f(\hat x_{i+1})}$ or $\Pi_{\hat f(\hat x_i)}$ intersects the interior of  $H^+_{\hat f(\hat x_{i+1})}-H^+_{\hat f(\hat x_{i-1})}$ (the planes $\Pi_{\hat f(\hat x_i)}$ are ``ordered'');
\item - $\Pi_{\hat f(\hat x_i)}\cap\Pi_{\hat f(\hat x_{i+1})}$ contains a geodesic $\hat d$ such that the nearest point retraction from $\hat d$  to $\hat f(\Hp^2)$ intersects the sub-arc of $\hat f(\hat k)$  joining $\hat f(\hat x_i)$ to $\hat f(\hat x_{i+1})$.\medskip
\end{description}

\indent
The integer $p-1$ (the number of components of $\hat k-\{\hat x_i\}$) is the {\it length} of the polygonal approximation. We will denote by $\theta(\Pi_{\hat f(\hat x_i)},\Pi_{\hat f(\hat x_{i+1})})$ the internal angle of $H^+_{\hat f(\hat x_{i})}-H^+_{\hat f(\hat x_{i+1})}$.\\
\indent
The existence of a polygonal approximation to any arc $\hat f(\hat k)$ intersecting at most once any leaf of $\hat f(\hat L)$ is proved in \cite{ceg}.\\
 
The {\em bending measure} $\int_{\hat k}\! d\hat\lambda$ along $\hat k$ is defined by $\int_{\hat k}\! d\hat\lambda=\inf_{\textgoth P}\sum_{i=1}^{p-1} \theta(\Pi_{\hat f(\hat x_{i})},\Pi_{\hat f(\hat x_{i+1})})$ where ${\textgoth P}$ runs over all polygonal approximations to $\hat f(\hat k)$.

It is shown in \cite[section 1.11]{sull} that this defines a transverse measure $\hat \lambda$ on the pleating locus of $\hat f$.\\

\indent A polygonal approximation ${\textgoth P}=\{(\hat x_i,\Pi_i)\}$ to a path $\hat f(\hat k)$ is an {\em $(\alpha,s)$-approximation} if
\begin{description}
\item - $\max_{1\leq i\leq p-1}\theta(\Pi_{i},\Pi_{i+1})<\alpha$ and
\item - $\max d_{\Hp^2}(\hat x_i,\hat x_{i+1})<s$.
\end{description}

\indent
The existence of a $(\delta,\eps)$-approximation for any $(\delta,\eps)$, and any arc $\hat f(\hat k)$ is shown in \cite{sull}. In the sequel we will need to have $(\delta,\eps)$-approximations with bounded length. The following lemma shows their existence.

\begin{lemme}   \label{beneu}
Let $\delta$ and $\eps$ be two positive numbers such that $0<\eps<\frac{\log 3}{2}$. Let $(\hat f,\Gamma,\rho)$ be a convex pleated surface and let $\hat k\subset\Hp^2$ be an arc that intersects the pleating locus $|\hat\lambda|$ transversely so that $\hat f(\hat k)$  intersects at most once any leaf of $\hat f(|\hat\lambda| )$. Then, there is a $(\delta,\eps)$-approximation to $\hat f(\hat k)$  with length at most $4(\frac{l(\hat k)}{\eps}+1)(\frac{\pi}{\delta}+1)=B(\eps,\delta,l(\hat k))$.
\end{lemme}

\begin{proof}
Consider the integer $p$ satisfying  $p-1\leq\frac{l(\hat k)}{\eps}\leq p$ and choose $p+1$ ordered points $\hat x_1,\ldots ,\hat x_{p+1}$ in $\hat k-|\hat\lambda|$ such that we have $\{\hat x_1,\hat x_{p+1}\}\subset\partial \hat k$ and $d(\hat x_i,\hat x_{i+1})\leq\eps$  for any $i\in\{1,\ldots,p+1\}$.  Choose also a support plane $\Pi_{\hat f(\hat x_i)}$ at $\hat f(\hat x_i)$ for each $i\in\{1,\ldots,p+1\}$. The first step of the proof will be to extend this family of support planes to obtain a polygonal approximation. There are three possible configuration for the positions of $\Pi_{\hat f(\hat x_i)}$ and $\Pi_{\hat f(\hat x_{i+1})}$. Let $\hat k_i$ be the sub-arc of $\hat k$ joining $\hat x_i$ to $\hat x_{i+1}$.\\

\indent - First configuration : $\Pi_{\hat f(\hat x_i)}$ intersects $\Pi_{\hat f(\hat x_{i+1})}$ and $\Pi_{\hat f(\hat x_i)}\cap\Pi_{\hat f(\hat x_{i+1})}$ contains a geodesic $\hat d$ such that the nearest point retraction from $\hat d$  to $\hat f(\Hp^2)$ intersects the arc $\hat f(\hat k_i)\subset \hat f(\hat k)$. In this configuration, $\{(\hat x_i,\Pi_{\hat f(\hat x_i)});(\hat x_{i+1},\Pi_{\hat f(\hat x_{i+1})})\}$ is already a polygonal approximation to $\hat f(\hat k_i)$.\\
\indent - Second configuration : $\Pi_{\hat f(\hat x_i)}$ does not intersect $\Pi_{\hat f(\hat x_{i+1})}$. Let $y$  be a point of $\hat k_i$ and let $\Pi_{\hat f(\hat y)}$ be a support plane at $\hat f(\hat y)$. The 3 half-spaces $H^-_{\hat f(\hat y)}$, $H^-_{\hat f(\hat x_i)}$ and $H^-_{\hat f(\hat x_{i+1})}$ intersect the ball $B(\hat f(\hat y),\eps)\subsetneq B(\hat f(\hat y),\frac{\log 3}{2})$. By \cite{gabai} (see also \cite{pont}) these 3 half-spaces are not disjoint, hence $\Pi_{\hat f(\hat y)}$ intersects either $\Pi_{\hat f(\hat x_i)}$ or $\Pi_{\hat f(\hat x_{i+1})}$. So any support plane at a point of $\hat f(\hat k_i)$ intersects either $\Pi_{\hat f(\hat x_i)}$ or $\Pi_{\hat f(\hat x_{i+1})}$. The arc $\hat f(\hat k_i)$ can be extended to an arc into the set of all the support planes at $\hat f(\hat k_i)$. Therefore there is a point $\hat y\subset int(\hat k_i)$ and a support plane  $\Pi_{\hat f(\hat y)}$ at $\hat f(\hat y)$ intersecting both $\Pi_{\hat f(\hat x_i)}$ and $\Pi_{\hat f(\hat x_{i+1})}$ (cf. \cite[\S 3.4]{kes}). Furthermore when we follow this arc $\Pi(t)$ in the set of support planes joining $\Pi_{\hat f(\hat x_i)}$ to $\Pi_{\hat f(\hat x_{i+1})}$, the nearest point retraction from $\Pi(t)\cap\Pi_{\hat f(\hat x_i)}$ ($\Pi(t)\cap\Pi_{\hat f(\hat x_i)}$ is a geodesic for $t>0$ small enough) to $\hat f(\Hp^2)$ moves from $\hat f(\hat x_i)$ in the direction of $\hat f(\hat x_{i+1})$ transversely to $\hat f(\hat k_i)$. It follows that the nearest point retraction from $\Pi_{\hat f(\hat x_i)}\cap\Pi_{\hat f(\hat y)}$  to $\hat f(\Hp^2)$ intersects the sub-arc of $\hat f(\hat k_i)$ joining $\hat f(\hat x_i)$ to $\hat f(\hat y)$. With a similar argument, we get that  the nearest point retraction from $\Pi_{\hat f(\hat y)}\cap\Pi_{\hat f(\hat x_{i+1})}$  to $\hat f(\Hp^2)$ intersects the sub-arc of $\hat f(\hat k_i)$ joining $\hat f(\hat y)$ to $\hat f(\hat x_{i+1})$.
Therefore $\{(\hat x_i,\Pi_{\hat f(\hat x_i)});(\hat y,\Pi_{\hat f(\hat y)});(\hat x_{i+1},\Pi_{\hat f(\hat x_{i+1})})\}$  is a polygonal approximation to $\hat f(\hat k_i)$.\\
\indent - Third configuration : $\Pi_{\hat f(\hat x_i)}$ intersects $\Pi_{\hat f(\hat x_{i+1})}$ but for any geodesic $\hat d\subset\Pi_{\hat f(\hat x_i)}\cap\Pi_{\hat f(\hat x_{i+1})}$, the nearest point retraction from $\hat d$  to $\hat f(\Hp^2)$ does not intersect $\hat f(\hat k_i)$. If any support plane at a point of $\hat f(\hat k_i)$ intersects either $\Pi_{\hat f(\hat x_i)}$ or $\Pi_{\hat f(\hat x_{i+1})}$ then we can find, as in the preceding case, a point $\hat y$  and a support plane $\Pi_{\hat f(\hat y)}$ such that \linebreak
 $\{(\hat x_i,\Pi_{\hat f(\hat x_i)});(\hat y,\Pi_{\hat f(\hat y)});(\hat x_{i+1},\Pi_{\hat f(\hat x_{i+1})})\}$  is a polygonal approximation to $\hat f(\hat k_i)$. Otherwise, there is a point $\hat y\in \hat k$  and a support plane  $\Pi_{\hat f(\hat y)}$ at $\hat f(\hat y)$ such that $\Pi_{\hat f(\hat y)}$  does not intersect $\Pi_{\hat f(\hat x_i)}$ nor $\Pi_{\hat f(\hat x_{i+1})}$. The planes $\Pi_{\hat f(\hat x_i)}$ and $\Pi_{\hat f(\hat y)}$ are in the second configuration, so there is a point $\hat z\subset \hat k_i$ between $\hat x_i$ and $\hat y$ such that $\{(\hat x_i,\Pi_{\hat f(\hat x_i)});(\hat z,\Pi_{\hat f(\hat z)});(\hat y,\Pi_{f(\hat y)})\}$ is a polygonal approximation to the sub-arc of $\hat k$ joining $\hat x_i$ to $\hat y$. Doing the same for $\Pi_{\hat f(\hat y)}$ and $\Pi_{\hat f(\hat x_{i+1})}$ we get a polygonal approximation to $\hat f(\hat k_i)$.\\

\indent 
Let us do the construction above for all the components of $\hat k-\{\hat x_1,\ldots,\hat x_{p}\}$. In each component of $\hat k-\{\hat x_1,\ldots,\hat x_{p}\}$, we have added at most $3$ points. So the resulting polygonal approximation has a length smaller than $4p$. Let us denote by\linebreak $\{(\hat x_i,\Pi_{\hat f(\hat x_i)})/ i=1,\ldots, 4p+1\}$ this polygonal approximation and let us denote by $\hat k_i\subset \hat k$ the geodesic arc joining $\hat x_i$ to $\hat x_{i+1}$. Consider $q\in\N$ satisfying $q-1\leq \frac{\theta(\Pi_{\hat f(\hat x_i)},\Pi_{\hat f(\hat x_{i+1})})}{\delta}\leq q$. Since $\theta(\Pi_{\hat f(\hat x_i)},\Pi_{\hat f(\hat x_{i+1})})\leq\pi$ we have $q\leq\frac{\pi}{\delta}+1$. We have already seen that we can extend an arc $\hat f(\hat k_i)$ to an arc in the set of all the support planes at $\hat f(\hat k_i)$. This implies that there are points $\hat y_j\in \hat k_i$, $1\leq j\leq q+1$, such that we have $\hat y_1=\hat x_i$, $\hat y_{q+1}=\hat x_{i+1}$ and $\theta(\Pi_{\hat f(\hat y_j)},\Pi_{\hat f(\hat y_{j+1})})\leq\delta$. Choosing such points for each arc $\hat k_i$, $1\leq i\leq 4p$, we get a $(\delta,\eps)$-approximation with length smaller than $4pq$. Since $p\leq\frac{l(\hat k)}{\eps}+1$ and $q\leq\frac{\pi}{\delta}+1$, this polygonal approximation satisfies the conclusion of Lemma \ref{beneu}.
\end{proof}

The following proposition of \cite{kes} gives an estimate of the error which is made when approximating the bending measure : 
\begin{prop}[KeS, Proposition 4.8]      \label{kesprop}
There is a universal constant $K$, and a function $s(\alpha)$, $0<s(\alpha)<1$, such that if ${\textgoth P}=\{(\hat x_i,\pi_i)\}$ is an $(\alpha,s(\alpha))$-approximation to a path $\hat f(\hat k)$, where $\alpha<\frac{\pi}{2}$, then we have
$$|\sum_{\textgoth P} \theta(\Pi_{i},\Pi_{i+1})-\int_k\! d\gamma|<K\alpha\, l(k)$$
\end{prop}

\indent Now we will use this proposition to prove the continuity of the bending measured geodesic lamination of a converging sequence of convex pleated surfaces.

\begin{lemme}   \label{plissee}
Let $(\hat f_n,\Gamma_n,\rho_n)$ be a sequence of convex pleated surfaces converging to a pleated surface $(\hat f_\infty,\Gamma_\infty,\rho_\infty)$ and let $\hat\lambda_n$ be the bending measured geodesic lamination of $\hat f_n$. The sequence $(\hat\lambda_n)$ converges for the weak$^*$ topology to a measured geodesic lamination $\hat\lambda_\infty$ and we have one of the following two situations :
\begin{description}
\item - $(\hat f_\infty,\Gamma_\infty,\rho_\infty)$ is a convex pleated surface and $\hat\lambda_\infty$ is its bending measured geodesic lamination;
\item - $(\hat f_\infty,\Gamma_\infty,\rho_\infty)$ is an even pleated surface, $|\hat\lambda_\infty|$ is the pleating locus of  $\hat f_\infty$ and $\hat\lambda_\infty$ is obtained by endowing each leaf of $|\hat\lambda_\infty|$ with a Dirac mass with a weight equal to $\pi$.
\end{description}
\end{lemme}
\begin{proof}
Notice that when $\Hp^3/\rho_\infty(\pi_1(M))$ is quasi-isometric to $\Hp^3/\rho_n(\pi_1(M))$ this result is a consequence of results of \cite{shear}.
By lemma \ref{confer}, $(\hat f_\infty,\Gamma_\infty,\rho_\infty)$ is an even or convex pleated surface.  We will show that any subsequence of $(\hat f_n,\Gamma_n,\rho_n)$ contains a subsequence satisfying the conclusions of lemma \ref{plissee}, the conclusion follows from this fact.\\
\indent
Let us choose a subsequence such that $(|\hat\lambda_n|)$ converge in the Hausdorff topology to a geodesic lamination $\hat L_\infty$. By \cite[\S 5.2]{ceg}, the pleating locus of $\hat f_\infty$ lies in $\hat L_\infty$. Let $\hat k\subset\Hp^2$ be an arc intersecting $\hat L_\infty$ transversely such that $\hat f_\infty(\hat k)$ intersects at most once any leaf of $\hat f_\infty(\hat L_\infty)$. Since $(\hat f_n(|\hat\lambda_n|))$ tends to $\hat f_\infty(\hat L_\infty)$ in the Hausdorff topology and since $(\hat f_n)$ converge to $\hat f_\infty$, for $n$ large enough $\hat f_n(\hat k)$ intersects at most once any leaf of $\hat f(|\hat \lambda_n|)$. Fix $\delta<\frac{\pi}{2}$ and $\eps<s(\delta)$ and choose for each $n$ a $(\delta,\eps)$-approximation $(\hat x_{i,n},\Pi_{\hat f_n(\hat x_{i,n})})$ to $\hat f_n(\hat k)$ whose length is the number $B(\delta,\eps,l(\hat k))$ appearing in Lemma \ref{beneu}. Extract a subsequence such that for any $i\leq B(\delta,\eps,l(\hat k))$, the sequence $(\hat x_{i,n},\Pi_{\hat f_n(\hat x_{i,n})})$ converges and denote by $(\hat x_{i,\infty},\Pi_{\hat f_\infty (\hat x_{i,\infty})})$ its limit. Since $\hat f_n$ converges to $\hat f_\infty$ and since $(\hat x_{i,n},\Pi_{\hat f_n(\hat x_{i,n})})$ converges to $(\hat x_{i,\infty},\Pi_{\hat f_\infty (\hat x_{i,\infty})})$, the nearest point retraction from $\Pi_{\hat f_n(\hat x_{i,n})}$ to $\hat f_n(\Hp^2)$ converges to the nearest point retraction from $\Pi_{\hat f_n(\hat x_{i,\infty})}$ to $\hat f_\infty(\Hp^2)$. It follows that $\{(\hat x_{i,\infty},\Pi_{\hat f_\infty (\hat x_{i,\infty})})\}$ satisfies all the requirements for being a $(\delta,\eps)$-approximation  to $\hat f_\infty(\hat k)$. Furthermore, the length of  $\{(\hat x_{i,\infty},\Pi_{\hat f_\infty (\hat x_{i,\infty})})\}$ is $B(\delta,\eps,l(\hat k))$.\\

\begin{claim}   \label{concon}
If $(\hat f_\infty,\Gamma_\infty,\rho_\infty)$ is a convex pleated surface, then $(\hat\lambda_n)$  converge to the bending measured geodesic lamination $\hat\lambda_\infty$ of $(\hat f_\infty,\Gamma_\infty,\rho_\infty)$.
\end{claim}

\begin{proof}
\indent Assume that $(\hat f_\infty,\Gamma_\infty,\rho_\infty)$ is a convex pleated surface and let us denote by $\hat\lambda_\infty$ its bending measured geodesic lamination. Let us recall that $\hat\lambda_n$ tends to $\hat\lambda_\infty$ if for any arc $\hat k\subset\Hp^2$ transverse to $|\hat\lambda_\infty|$, $(\int_{\hat k}\! d\hat\lambda_n)$ tends to $\int_{\hat k}\! d\hat\lambda_\infty$.\\
\indent Let $\hat k\subset\Hp^2$ be an arc transverse to $|\hat\lambda_\infty|$ and let $\{(\hat x_{i,n},\Pi_{\hat f_n(\hat x_{i,n})})\}$ be the $(\delta,\eps)$-approximations defined above.
 For any $i\leq p$, $\Pi_{\hat f_n(\hat x_{i,n})}$ tends to $\Pi_{\hat f_\infty (\hat x_{i,\infty})}$. Therefore, for $n$ large enough, we have $\sum_{i=1}^{p-1} |\theta(\Pi_{\hat f_\infty (\hat x_{i,\infty})},\Pi_{\hat f_\infty (\hat x_{i+1,\infty})})-\theta(\Pi_{\hat f_n(\hat x_{i,n})},\Pi_{\hat f_n(\hat x_{i+1,n})})|\leq \delta$. It follows from Proposition \ref{kesprop} that we have $|\sum_{i=1}^{p-1} \theta(\Pi_{\hat f_n(\hat x_{i,n})},\Pi_{\hat f_n(\hat x_{i+1,n})})-\int_{\hat k}\! d\hat\lambda_n|<K \delta l(\hat k)$ for any $n\in\overline{\N}=\N\cup\{\infty\}$. Hence we have $|\int_{\hat k}\! d\hat\lambda_\infty-\int_{\hat k}\! d\hat\lambda_n|<2K \delta l(\hat k)+\delta$. Letting $\delta$ tend to $0$ yields $\int_{\hat k}\! d\hat\lambda_n\rightarrow\int_{\hat k}\! d\hat\lambda_\infty$. The same is true for any arc $\hat k$ transverse to $|\hat\lambda_\infty|$, hence $\hat\lambda_n$ converges to $\hat\lambda_\infty$.
\end{proof}

\indent Next we will consider the case where $(\hat f_\infty,\Gamma_\infty,\rho_\infty)$ is an even pleated surface.  Let us first show that all the leaves of the pleating locus $|\hat\lambda_\infty|$ of $(\hat f_\infty,\Gamma_\infty,\rho_\infty)$ are isolated leaves, namely that their projections to $\Hp^2/\Gamma_\infty$ are isolated leaves.

\begin{claim}   \label{isol}
When $(\hat f_\infty,\Gamma_\infty,\rho_\infty)$ is an even pleated surface, the pleating locus $|\hat\lambda_\infty|$ of $\hat f_\infty$ contains only isolated leaves.
\end{claim}
\begin{proof}
Choose two distinct successive points $\hat x_{i,\infty}$ and $\hat x_{i+1,\infty}$ and denote by $\hat k_i$ the sub-arc of $\hat k$ joining $\hat x_{i,\infty}$ to $\hat x_{i+1,\infty}$. Let $\hat y$ be a point of $int(\hat k_i)$ and let $(\Pi_{\hat f_n(\hat y)})$ be a sequence of support planes at $\hat f_n(\hat y)$ converging to a support plane $\Pi_{\hat f_\infty(\hat y)}$ at $\hat f_\infty(\hat y)$. For $n$ large enough, $\hat y$ lies in the sub-arc of $\hat k$ joining $\hat x_{i,n}$ to $\hat x_{i+1,n}$. Since $\{(\hat x_{i,n},\Pi_{\hat f_n(\hat x_{i,n})})\}$ is a $(\delta,\eps)$-approximation, $\Pi_{\hat f_n(\hat y)}$ intersects both $\Pi_{\hat f_n(\hat x_{i,n})}$ and $\Pi_{\hat f_n(\hat x_{i+1,n})}$ and we have $\theta(\Pi_{\hat f_n(\hat x_{i,n})},\Pi_{\hat f_n(\hat y)})~\leq~\delta$ and $\theta(\Pi_{\hat f_n(\hat y)},\Pi_{\hat f_n(\hat x_{i+1,n})})~\leq~\delta$. Letting $n$ tend to $\infty$ yields  $\theta(\Pi_{\hat f_\infty(\hat y)},\Pi_{\hat f_\infty(\hat x_{i,\infty})})\leq\delta$ and $\theta(\Pi_{\hat f_\infty(\hat y)},\Pi_{\hat f_\infty(\hat x_{i+1,\infty})})\leq\delta$. Since $(\hat f_\infty,\Gamma_\infty,\rho_\infty)$ is an even pleated surface, the dihedral angle between two adjacent support planes lies in $\{0,\pi\}$. Hence we have
 $\theta(\Pi_{\hat f_\infty(\hat y)},\Pi_{\hat f_\infty(\hat x_{i,\infty})})=0$ and $\theta(\Pi_{\hat f_\infty(\hat y)},\Pi_{\hat f_\infty(\hat x_{i+1,\infty})})=0$.  It follows that for any arc
 $\hat\kappa\subset int (\hat k_i)$, $\int_{\hat\kappa}\! d\hat\lambda$ tends to $0$. From the proof of  \cite[Prop 27]{shear} we deduce that $int (\hat k_i)$ does not intersect the pleating locus of $\hat f_\infty$.\\
\indent
So we have shown that the intersection between $\hat k$ and the pleating locus of $\hat f_\infty$ lies in $\{\hat x_i|\; i=1\ldots p\}$, in particular this intersection has a finite cardinal.
\end{proof}

Endowing each leaf of $|\hat\lambda_\infty|$ with a Dirac measure whose weight is equal to $\pi$ yields a measured geodesic lamination $\hat\lambda_\infty$.\\
\indent
Let $r$ be a lower bound for the set $\{d(\hat x_1,\hat x_2)/ \hat x_1\mbox{ and } \hat x_2$ are two different points of $\hat k\cap|\hat\lambda_\infty| \}$. By the proof of Claim \ref{isol}, we can choose $r>0$. If $\{(\hat x_{i,\infty},\Pi_{\hat f_\infty (\hat x_{i,\infty})})\}$ is a $(\delta,\eps)$-approximation satisfying $\eps<r$, we have $ \sum_{i=1}^{q-1} \theta(\Pi_{\hat f_\infty (\hat x_{i,\infty})},\Pi_{\hat f_\infty (\hat x_{i+1,\infty})})=\int_{\hat k}\! d\hat\lambda_\infty$.\\
\indent In the proof of Claim  \ref{concon}, we can use this equality instead of  the inequality of proposition \ref{kesprop}. Thus we get the following claim which concludes the proof of Lemma \ref{plissee} : 

\begin{claim}   \label{verg}
The sequence $(\hat\lambda_n)$ tends to $\hat\lambda_\infty$ in ${\cal ML}(\Hp^2)$.
\end{claim}

Let us notice that in the case where $(\hat f_\infty,\Gamma_\infty,\rho_\infty)$ is even, we did not show that the projection of  $|\hat\lambda_\infty|$ to $\Hp/\Gamma_\infty$  is compact. We will see further in the text that this is actually true.
\end{proof}

\indent
We will now improve the description of $(\hat f_\infty,\Gamma_\infty,\rho_\infty)$ when it is an even pleated surface.

\begin{lemme}   \label{plat}
Let $(\hat f_n,\Gamma_n,\rho_n)$ be a sequence of convex pleated surfaces converging to an even pleated surface $(\hat f_\infty,\Gamma_\infty,\rho_\infty)$. There is a surface $S$ with geodesic boundary, such that $\Hp^2/\Gamma_\infty$ is the double of $S$ and such that the pleating locus $|\hat\lambda_\infty|$ of $\hat f_\infty$ project to $\partial S\subset \Hp/\Gamma_\infty$. Furthermore $\partial S$ is compact.
\end{lemme}

\begin{proof}
First the following claim will define $S$.

\begin{claim}   \label{mesu}
All the connected components of $\Hp^2-|\hat\lambda_\infty|$ have the same image under $\hat f_\infty$.
\end{claim}

\begin{proof}
Recall that since $\hat f_\infty$ is a pleated surface, the restriction of $\hat f_\infty$ to each connected component of $\Hp^2-|\hat\lambda_\infty|$ is one to one.\\
\indent Assume that Claim \ref{mesu} is not true. There exist two connected components $\hat P$ and $\hat P'$ of $\Hp^2-|\hat\lambda_\infty|$ whose closures intersect and whose images under $\hat f_\infty$ are different. There is a point $\hat y$ lying in the boundary of the closure of $\hat P$ such that $f_\infty(\hat y)$ lies in $\hat f_\infty(\hat P')$  or a point $\hat y$ lying in the boundary of the closure of $\hat P'$ such that $\hat f_\infty(\hat y)$ lies in $f_\infty(\hat P)$. This two cases are similar and we will only deal with the first one. Let us denote by $\hat l$ the leaf of $|\hat\lambda_\infty|$ containing $\hat y$. By the proof of Claim \ref{isol}, $\hat l$ is an isolated leaf. Let $\hat k\subset\Hp^2$ be a geodesic arc intersecting $\hat l$ transversally so that $\hat k\cap |\hat\lambda_\infty|=\{y\}$. For $n\in{\overline \N}$, let $(\hat x_{i,n},\Pi_{\hat f_n (\hat x_{i,n})})$ be the  $(\delta,\eps)$-approximation to $\hat f_n(\hat k)$ constructed at the beginning of the proof of Lemma \ref{plissee} and  let us denote by $p=B(\delta,\eps,l(\hat k))$ its length. By Claim \ref{verg}, the bending measure on $\hat k$ tends to a Dirac measure whose weight is equal to $\pi$ and whose support is $\{y\}$. It follows that there is $i$ such that $(\hat x_{i,n})$ and $(\hat x_{i+1,n})$ converge to $\hat y$ and such that we have the inequality $\frac{\pi}{p}\leq\theta(\Pi_{\hat f_n (\hat x_{i,n})},\Pi_{f_n (\hat x_{i+1,n})})\leq\delta$ for $n$ large enough.\\
\indent
Let $\hat z$ be the point of $\hat P'$ such that we have $\hat f_\infty(\hat z)=\hat f_\infty(\hat y)$. For any $n$, $\hat z$ lies either in a connected component $\hat P'_n$ of $\Hp^2-|\hat\lambda_n|$ or in the closures of two connected components of $\Hp^2-|\hat\lambda_n|$, in this second case we will denote by $\hat P'_n$ the interior of the union of these two closures. Let us extract a subsequence such that $\hat P'_n$ converges to an open subsurface $\hat P'_\infty$ of $\hat P'$ (for example a subsequence such that $|\hat\lambda_n|$ converges in the Hausdorff topology). Since $\hat f_n$ converge to $\hat f_\infty$ on any compact set, the sequence $\hat f_n(P'_n)$ converges to\linebreak $\hat f_\infty(P'_\infty)\subset\hat f_\infty(P')$.\\
\indent The point $\hat f_\infty(\hat y)$ lies in  $\hat f_\infty(\hat P'_\infty)$, in $\Pi_{\hat f_\infty (\hat  x_{i,\infty})}$ and in $\Pi_{\hat f_\infty (\hat x_{i+1,\infty})}$ and we have the inequality $\frac{\pi}{p}\leq\theta(\Pi_{\hat f_n (\hat x_{i,n})},\Pi_{\hat f_n (\hat x_{i+1,n})})\leq\delta$. Therefore, for $n$ large enough, one plane among $\Pi_{\hat f_n (\hat x_{i,n})}$ and $\Pi_{\hat f_n (\hat x_{i+1,n})}$ intersects $\hat f_n(\hat P'_n)$ transversely. This contradicts the convexity of $(\hat f_n,\Gamma_n,\rho_n)$ and concludes this proof.
\end{proof}

Next we will show that the quotient of $\Hp^2-|\hat\lambda_\infty|$ by $\Gamma_\infty$ has two connected components.\\
\indent
Let us extract a subsequence such that $(|\hat\lambda_n|)$ converges in the Hausdorff topology to a geodesic lamination $\hat L_\infty$ and let $\hat x,\hat y\subset \Hp^2-\hat L_\infty$ be two points whose images $\hat f_\infty(\hat x)$ and $\hat f_\infty(\hat y)$ coincide. Denote by $\hat P_{\hat x}$ and $\hat P_{\hat y}$ the connected components of $\Hp^2-|\hat\lambda_\infty|$ containing $\hat x$ and $\hat y$ and assume that $H^+_{\hat f_\infty(\hat P_{\hat x})}$is equal to $H^+_{\hat f_\infty(\hat P_{\hat y})}$. We will show the existence of an element $g_\infty$ of $\Gamma_\infty$ such that $\hat x=g_\infty\hat y$.\\
\indent
The half-spaces $H^+_{\hat f_n(\hat P_{\hat x})}$ and $H^+_{\hat f_n(\hat P_{\hat y})}$ converge to $H^+_{\hat f_\infty(\hat P_{\hat x})}$ and we have $\hat f_n(\hat P_{\hat y})\subset H^+_{\hat f_n(\hat P_{\hat x})}$ and
 $\hat f_n(\hat P_{\hat x})\subset H^+_{\hat f_n(\hat P_{\hat y})}$. It follows that the planes $\Pi_{\hat f_n(\hat x)}$ and $\Pi_{\hat f_n(\hat y)}$ intersect each other for $n$ large enough and that the sequence $\theta(\Pi_{\hat f_n(\hat x)},\Pi_{\hat f_n(\hat y)})$ tends to $0$. This implies  that the distance $d_{(\Pi_{\hat f_n(\hat x)}\cup\Pi_{\hat f_n(\hat y)})}(\hat f_n(\hat x),\hat f_n(\hat y))$ measured on $\Pi_{\hat f_n(\hat x)}\cup\Pi_{\hat f_n(\hat y)}$ tends to $0$. Since we have\linebreak $d_{(\Pi_{\hat f_n(\hat x)}\cup\Pi_{\hat f_n(\hat y)})}(\hat f_n(\hat x),\hat f_n(\hat y))\geq d_{\hat f_n(\Hp^2)}(\hat f_n(\hat x),\hat f_n(\hat y))$ and since $\hat f_\infty$ is a local homeomorphism on the complementary regions of $|\hat\lambda_\infty|$, there is, for $n$ large enough, an isometry $a_n\in\Gamma_n$ such that $d_{\Hp^2}(a_n \hat y,\hat x)\longrightarrow 0$. The distance $d_{\Hp^2}(\hat y,a_n\hat y)$ is bounded, hence, up to extracting a subsequence, $(a_n)$ converges to some $a_\infty\in\Gamma_\infty$. We have then $a_\infty\hat y=\hat x$.\\
\indent
This implies that any point of $\hat f_n (\Hp^2-|\hat\lambda_\infty|)$ has at most two preimages in $\Hp^2/\Gamma_\infty$. Combining this with Claim \ref{mesu}, we get the existence of a surface $S$  such that $\Hp^2/\Gamma_\infty$ is the double of $S$ and that the pleating locus $|\hat\lambda_\infty|$ of $\hat f_\infty$ projects to $\partial S\subset \Hp^2/\Gamma_\infty$.\\
\indent
It remains to show that $\partial S$ is compact. Assume the contrary, then $\partial S$ contains two asymptotic half geodesics. Let $\hat k\subset S$ be a geodesic arc joining these two half geodesic. The double of $\hat k$ is a simple closed curve $c\subset\Hp^2/\Gamma_\infty$ bounding a cusp of $\Hp^2/\Gamma_\infty$. The image of $c$ under $f_\infty$ is $f_\infty(\hat k)$ covered twice. It follows that the curve $f_\infty(c)$ is homotopic to a point in $\Hp^3/\rho_\infty(\Gamma_\infty)$. This contradicts the assumption that the parabolic elements of $\Gamma_\infty$ are mapped to parabolic isometries by $\rho_\infty$.
\end{proof}

\indent
We will conclude this section with two lemmas which will be used in the next sections. They are proved in \cite{espoir} (see also \cite{petitbain}). We will call them slight bending Lemmas.

\begin{lemme}   \label{lemap}
Let $(\hat f,\Gamma,\rho)$ be a convex pleated surface, let $\hat x_1$ and $\hat x_2\in\Hp^2$, let $\hat c_1\subset\Hp^2$ be the geodesic segment joining $\hat x_1$ to $\hat x_2$ and let $\tilde c_2\subset\Hp^3$ be the geodesic segment joining $\hat f(\hat x_1)$ to $\hat f(\hat x_2)$. If there exists $\eps<\frac{\pi}{2}$ such that the bending measure of $\hat f(\hat c_1)$ is smaller than $\eps$, then $\exists C_{\eps}$ such that
 $l(\hat{c}_1)\leq C_{\eps} l(\tilde{c}_2)$. Furthermore $\lim_{\eps\longrightarrow 0} C_{\eps}= 1$ and the sum of the exterior angles that $\hat f(\hat c_1)$ and $\tilde c_2$ make at their vertices is smaller than $\eps$.
\end{lemme}

\begin{lemme}[Slightly bent curves are quasi-geodesics] \label{approx}
Let $(\hat f,\Gamma,r)$ be a convex pleated surface, let $\hat\lambda$ be its bending measured geodesic lamination, let $c$ be a simple closed geodesic of $\Hp^2/\Gamma$ and let $c^*$ be the geodesic of $\Hp^3/r(\Gamma_n)$ in the homotopy class of $f(c)$. For any $\eps<\frac{\pi}{2}$ there exists $C_{\eps}$ and $A_{\eps}$ such that if $i(c,\lambda)\leq\eps$ then $l(c)\leq C_\eps(l(c^*)+A_\eps)$. Moreover $\lim_{\eps\longrightarrow 0} C_\eps=1$, and $\lim_{\eps\longrightarrow 0} A_\eps=0$. 
\end{lemme}

\section{The continuity of $b_{\cal R}$}
Next we will use the results of the previous section to show the continuity of $b_{\cal R}$. But first let us precise the definition of the bending measured lamination of a geometrically finite metric $\sigma\in{\cal GF}(M)$. Let $\rho:\pi_1(M)\rightarrow Isom(\Hp^3)$ be a representation associated to $\sigma$ and let $N(\rho)$ be the Nielsen core of $\rho$. There are a multi-curve $\lambda^{(p)}$ and a natural (relative to $\sigma$) homeomorphism $f:\partial_{\chi<0} M-\lambda^{(p)}\rightarrow \partial N(\rho)$. This homeomorphism $f$ is constructed by using the retraction map from the domain of discontinuity to the boundary of the convex core of the limit set. It is well-defined up to isotopy. The surface $\partial N(\rho)$ is the image of a convex pleated surface and therefore, it has a bending measured lamination. Let us denote by $\lambda\in{\cal ML}(\partial M)$ the image under $f^{-1}$ of this bending measured lamination. Adding the leaves of $\lambda^{(p)}$ endowed with Dirac masses with weights equal to $\pi$ we get the bending measured geodesic lamination of $\sigma$. This measured lamination $\lambda$ does not depend on the choice of $\rho$ (among the representations associated to $\sigma$).\\
\indent
Given a simple closed geodesic $c$ and a hyperbolic metric $s$ on $\partial_{\chi<0} M$, let $l_s(c)$ be the length of the corresponding $s-geodesic$. If $\gamma$ is a weighted simple closed geodesic with a weight $w(\gamma)$ we define $l_s(\gamma)$ by $l_s(\gamma)=w(\gamma) l_s(|\gamma|)$. This function $l_s$ extends continuously to a function $l_s:{\cal ML}(\partial M)\rightarrow\R$.\\
\indent
Given a simple closed curve $c\in\partial_{\chi<0} M$, and a hyperbolic metric $\sigma$ on $int(M)$, we denote by $c^*$ the closed $\sigma$-geodesic in the free homotopy class of $c$ and by $l_\sigma(c^*)$ its $\sigma$-length.\\
\indent
The following proposition shows the continuity of the bending map $b_{\cal R}$. It is essentially a rephrasing of the statement of Theorem $1$. Thus Theorem $1$ follows from Proposition \ref{con}.

\begin{prop}    \label{con}
Let $(\sigma_n)$ be a sequence of geometrically finite metrics on the interior of $M$ converging to a non fuchsian geometrically finite metric $\sigma_{\infty}$ and let $\lambda_n$  be the bending measured geodesic lamination of $\sigma_n$. The sequence $\dot\lambda_n$ converges to $\dot\lambda\in{\cal P}(M)/{\cal R}$ and the bending measured geodesic lamination of $\sigma_\infty$ is the representative $\lambda'$ of $\dot\lambda$ that lies in ${\cal P}(M)$.
\end{prop}

\begin{proof}
We are going to show that any subsequence of $(\sigma_n)$ contains a subsequence satisfying the conclusion of the proposition.
Let us begin by a result about the curves whose length tends to $0$ when $n$ tends to $\infty$.
\begin{lemme}   \label{pandit}
Let $(\sigma_n)$ be a sequence of geometrically finite metrics on the interior of $M$ converging algebraically and let $(s_n)$ be the sequence of metrics induced on $\partial M$ by a homeomorphism $h_n : M\rightarrow N(\rho_n)^{ep}$ associated to $\sigma_n$. Then, there is $K>0$ such that if $c\subset\partial M$ is a simple closed curves that bounds an essential disc, we have $l_{s_n}(c)\geq K$.
\end{lemme}
\begin{proof}
Let $D\subset M$ be an essential disc and let $e\subset\partial M$ be a simple closed curve which does not bound an essential disc. We will say that $e$ intersects $\partial D$ {\em essentially} if the ends of any lift of $e$ to $\tilde M$ are separated by a lift of $\partial D$. As in the introduction ${\cal P}(M)$ consists of the measured geodesic laminations satisfying conditions $a)$, $b)$ and $c)$. 

\begin{claim}   \label{ess}
Let $\gamma\subset{\cal P}(M)$ be a weighted multi-curve and let $D\subset M$ be an essential disc; then at least one leaf of $\gamma$ intersects  $\partial D$ essentially.
\end{claim}
\begin{proof}
 Let $\partial\tilde D$ be a lift of $\partial D$ to $\partial\tilde M$. This lift separates $\partial\overline{\tilde M}\approx S^2$ into two discs $C$ and $C'$. Let us denote by $I_i,1\leq i\leq p$ (resp. $I'_i,1\leq i\leq p'$)  the connected components of $p^{-1}(\gamma)-\partial\tilde D$ lying in $C$ (resp. $C'$) whose endpoints lie in $\partial\tilde D$. These arcs $I_i$ (resp. $I'_i$) cut $C$ (resp. $C'$) in $p+1$ discs  $C_i$, $0\leq i\leq p$ (resp. in $p'+1$ discs $C'_i$, $0\leq i\leq p'$). The curves $\partial C_i$ and $\partial C'_i\subset\partial\tilde M$ bound essential discs in $\tilde M$. It follows from conditions $a)$ and $c)$ that for any $i$, we have $\sharp\{\partial C_i\cap p^{-1}(\gamma)\}\geq 3$  and  $\sharp\{\partial C'_i\cap p^{-1}(\gamma)\}\geq 3$. If no leaf of $\gamma$ intersects $\partial D$ essentially, then, for any $i$, each point of $\partial C_i\cap\partial\tilde D\cap p^{-1}(\gamma)$ is the endpoint of an $I'_j$ and each point of  $\partial C'_i\cap\partial\tilde D\cap p^{-1}(\gamma)$ is the endpoint of an $I_j$. This implies that $3(p+1)\leq 2p'$ and that $3(p'+1)\leq 2p$. We get $p+p'\leq -6$. This contradiction shows that at least one leaf of $\gamma$ intersects  $\partial D$ essentially.
\end{proof}

\indent
Let $( c_n)$ be a sequence of simple closed curves such that each $ c_n$ bounds an essential disc and that $l_{s_n}( c_n)\longrightarrow 0$. For large $n$, $ c_n$ is the core of a wide Margulis tube $T_n\subset\partial M$. Let  $c_{1,n}$ and $c_{2,n}$ be the boundary components of $T_n$ and let us choose $T_n$ so that we have\linebreak $l_{s_n}(c_{1,n})=l_{s_n}(c_{2,n})\longrightarrow 0$ and $d_{s_n}(c_{i,n},c_n)\longrightarrow\infty$. Let $e$ be a simple closed curve supporting a transverse measure $\gamma\in{\cal P}(M)$. By Claim \ref{ess}, $e$ intersects all the $ c_n$ essentially. For each $n$, let $e_n^*$ be the $\sigma_n$-geodesic in the homotopy class of $e$, let $\tilde e^*_n$, $\tilde c_n$ and $\tilde T_n$ be lifts of $e^*_n$, $ c_n$, $T_n$ respectively. Let $\tilde c_{1,n}$ and $\tilde c_{2,n}$ be the components of $\partial\tilde T_n$. Let $\Pi$ be a geodesic plane containing $\tilde e^*_n$. Let $\tilde k_n$ and $\tilde k'_n$ be the two connected components of $\tilde T_n\cap\Pi$ and let $\tilde\kappa_n$ and $\tilde\kappa'_n$ be the geodesic segment of $\Pi$ joining the endpoints of $\tilde k_n$ and $\tilde k'_n$, see figure \ref{meridien}.
\begin{figure}[hbtp]
\psfrag{a}{$\tilde k_n$}
\psfrag{b}{$\tilde\kappa_n$}
\psfrag{c}{$\tilde e^*$}
\psfrag{d}{$\tilde\kappa'_n$}
\psfrag{e}{$\tilde k'_n$}
\psfrag{f}{$\tilde d_n$}
\psfrag{i}{$\tilde d'_n$}
\psfrag{h}{$\alpha_2$}
\psfrag{g}{$\alpha_1$}
\centerline{\includegraphics{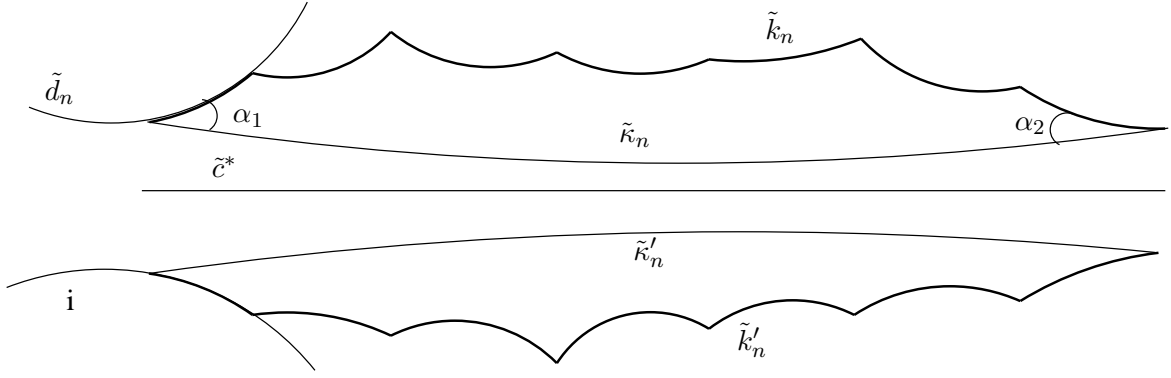}}
\caption{The section $\Pi\cap C(\rho_n)$}
\label{meridien}
\end{figure}
Assume that $\partial \tilde k_n$ and $\partial \tilde k'_n$ does not lie in the preimage of $\lambda_n$. Let $\tilde d_n$ and $\tilde d'_n$ be the intersections of $\Pi$ and of support planes at $\tilde k_n\cap\tilde c_{1,n}$ and at $\tilde k'_n\cap\tilde c_{1,n}$ respectively. The section $C(\rho_n)\cap\Pi$ lies between $\tilde d_n$ and $\tilde d'_n$. Since $\tilde e^*_n$ lies in $C(\rho_n)$, $\tilde e^*_n$ does not intersect $\tilde d_n$ nor $\tilde d'_n$. The arcs $\tilde\kappa_n$ and $\tilde\kappa'_n$ are very close, this implies that the angle $\alpha_1$ between $\tilde d_n$ and $\tilde\kappa_n$ is small. The same considerations apply for the angle $\alpha_2$ on the other vertex of $\tilde\kappa_n$. Let $\alpha<\frac{\pi}{2}$ be an upper bound for $\alpha_1$ and $\alpha_2$ and set $r_n=\max\{d(\tilde x,\tilde\kappa_n)|\tilde  x\in \tilde k_n\}\leq\max\{d(\tilde x,\tilde e^*_n)/ \tilde x\in \tilde k_n\}$. Any point $\tilde x$ lying in $\tilde k_n$ lies in a simple closed curve $\tilde c_{x,n}\subset\tilde T_n$ homotopic to $\tilde c_n$ such that $l_{s_n}(c_{x,n})\leq l_{s_n}(c_{1,n})$. This curve $\tilde c_{x,n}$ bounds a disc with diameter less than $l_{s_n}(c_{x,n})$ which intersects $\tilde e^*_n$. It follows that $r_n\leq l_{s_n}(c_{1,n})\longrightarrow 0$. A computation using Fermi coordinates (compare with \cite[Lemme A.2]{espoir}) yields $l_{\sigma_n}(\tilde k_n)\leq\sqrt{1+\tan^2\alpha}\,(\ch r_n)^2 l_{\sigma_n}(\tilde\kappa_n)$. Since the width of $T_n$ tends to infinity, we have $l_{\sigma_n}(\tilde k_n)\longrightarrow\infty$. So we get $l_{\sigma_n}(\tilde\kappa_n)\longrightarrow\infty$. Since any point of $\kappa_n$ is close to the geodesic $e^*_n$, we have $l_{\sigma_n}(e^*_n)\longrightarrow\infty$. But this contradicts the algebraic convergence of $(\sigma_n)$.
\end{proof}

\indent
Since $(\sigma_n)$ converges to $\sigma_{\infty}$, we may choose representations $\rho_n:\pi_1(M)\rightarrow Isom(\Hp^3)$ associated to $\sigma_n$ such that  $(\rho_n)$ converges algebraically to a representation $\rho_{\infty}$ associated to $\sigma_{\infty}$ and that $(\rho_n(\pi_1(M)))$ converges geometrically to $\rho_\infty(\pi_1(M))$. For $n\in\overline{\N}$, let $\lambda_n$ be the bending measured geodesic lamination of $\sigma_n$. Denote by $\lambda_n^{(p)}$ the union of the leaves of $\lambda_n$ which have a weight equal to $\pi$. Let $M_n=\Hp^3/\rho_n(\pi_1(M))$ and let $h_n : M-\lambda_n^{(p)}\rightarrow N(\rho_n)$ be a homeomorphism associated to $\sigma_n$.\\
\indent
Let $x_\infty$ be a point in $\partial N(\rho_{\infty})$, let $\partial_{x_\infty} N(\rho_\infty)$ be the connected component of $\partial N(\rho_\infty)$ containing $x_\infty$ and let $S$ be a connected component of $\partial M-\lambda_\infty^{(p)}$ satisfying\linebreak $h_\infty (S)=\partial_{x_\infty} N(\rho_\infty)$. Let $\tilde x_\infty\in C(\rho_\infty)$ be a lift of $x_\infty$ and let $\partial_{\tilde x_\infty} C(\rho_\infty)$ be the connected component of $\partial C(\rho_\infty)$ containing $\tilde x_\infty$. By \cite{taylor}, the convex hull $C(\rho_{\infty})$ of $L_{\rho_{\infty}}$ is the limit of $(C(\rho_n))$ with respect to the Hausdorff topology. Hence there exists $\tilde x_n\subset\partial C(\rho_n)$ such that $(\tilde x_n)$ converges to $\tilde x_\infty$. Let us denote by $\partial_{\tilde x_n} C(\rho_n)$ the connected component of $\partial C(\rho_n)$ containing $\tilde x_n$ and by $\partial_{\tilde x_n} N(\rho_n)$ its projection to $N(\rho_n)$. There is a connected component $S_n$ of $\partial M-\lambda_n^{(p)}$ and a pleated map $g_n:S_n\rightarrow M$ onto $\partial_{\tilde x_n} N(\rho_n)$. Let $(\hat g_n,\Gamma_n,r_n)$ be the pleated surface that lifts $g_n$, namely the map $\hat g_n:\Hp^2\rightarrow\partial_{\tilde x_n} C(\rho_n)$ is onto, $r_n:\Gamma_n\rightarrow\rho_n(\pi_1(M))$ is onto the stabiliser of $\partial_{\tilde x_n} C(\rho_n)$ and the quotient map $\Hp^2/\ker(r_n)\rightarrow\partial_{\tilde x_n} C(\rho_n)$ is a homeomorphism. By \cite[Theorem 5.2.2]{ceg} and Lemma \ref{pandit} a subsequence of the sequence of pleated surfaces $(\hat g_n,\Gamma_n,r_n)$ converges to a pleated surface $(\hat g_{\infty},\Gamma_\infty,r_\infty)$.\\
\indent
By Lemma \ref{confer}, $(\hat g_{\infty},\Gamma_\infty,r_\infty)$ is a convex or even pleated surface. Since $\tilde x_\infty\in \hat g_\infty(\Hp^2)$, $\hat g_\infty(\Hp^2)\subset\partial_{\tilde x_\infty} C(\rho_\infty)$ and $r_\infty(\Gamma_\infty)\subset\rho_\infty(\pi_1(M))$ is a subgroup of the stabiliser of $\partial_{\tilde x_\infty} C(\rho_\infty)$.

\begin{lemme}   \label{momeo}
The quotient map $g_\infty:\Hp^2/\Gamma_\infty\rightarrow \partial_{x_\infty} N(\rho_\infty)$ is a homeomorphism.
\end{lemme}
\begin{proof}
Set $C_{\hat g_n}=\bigcap\{H_{P_i}^{\epsilon}/P_i$ is a connected component of $\Hp^2-\hat L\}$. For any $n\in\N$, we have $C(\rho_n)\subset  C_{\hat g_n}$. Since $\sigma_n$ converges to $\sigma_\infty$, it follows from \cite{taylor} that $C(\rho_n)$ converges to $C(\rho_\infty)$ in the Hausdorff topology. So we get $C(\rho_\infty)\subset C_{\hat g_\infty}$. Since $\rho_\infty$ is not fuchsian, we have $int(C(\rho_\infty))\neq\emptyset$. Therefore $int(C_{\hat g_\infty})\neq\emptyset$ and $(\hat g_{\infty},\Gamma_\infty,r_\infty)$ is a convex pleated surface. By \cite{sull}, this implies that $g_\infty$ is a covering map.\\
\indent
Let $\hat\gamma_n$ be the bending geodesic lamination of $\hat g_n$. Extract a subsequence such that $(|\hat\gamma_n|)$ converges to some geodesic lamination $\hat L_\infty$ in the Hausdorff topology.\\
\indent
Assume that $g_\infty$ is not a homeomorphism. There are two points $z\neq y\in(\Hp^2-\hat L_\infty)/r_\infty(\Gamma_\infty)$ such that $g_\infty (z)=g_\infty(y)$. If we lift the situation to $\Hp^2$, we get $\hat z,\hat y\subset\Hp^2-\hat L_\infty$ and $\rho_\infty (a)\in \rho_\infty(\pi_1(M))$ such that $\hat z\not\in\Gamma_\infty\hat y$ and that $\hat g_\infty (\hat z)=\rho_\infty (a)\circ\hat g_\infty (\hat y)$. The sequences $\tilde z_n=\hat g_n(\hat z)$ and $\tilde y_n=\rho_n(a)\circ\hat g_n(\hat y)$ converge simultaneously to $\hat g_\infty (\hat z)$. Let $\Pi_{\tilde z_n }$ and $\Pi_{\tilde y_n}$ be support planes of $C_n(\rho_n)$ at $\tilde z_n $ and $\tilde y_n$. The half-spaces $(H^+_{\Pi_{\tilde z_n }})$ and $(H^+_{\Pi_{\tilde y_n}})$ tend to half-spaces $H^+_{\hat z}$ and $H^+_{\hat y}$ respectively.  These two half-spaces $H^+_{\hat z}$ and $H^+_{\hat y}$ are bounded by the support plane $\Pi_{\hat g_\infty (\hat z)}$ at $\hat g_\infty (\hat z)$. Thus we have either $H^+_{\tilde z}\cap H^+_{\tilde y}=\Pi_{\hat g_\infty (\hat z)}$, or $H^+_{\tilde z}= H^+_{\tilde y}$.\\
\indent
Since $C(\rho_\infty)\subset C_{\hat g_\infty}\subset H^+_{\tilde z_n}\cap H^+_{\tilde y_n}$ and since $\rho_\infty(\pi_1(M))$ is not fuchsian, we have\linebreak $H^+_{\tilde z}\cap H^+_{\tilde y}\neq\Pi_{\hat g_\infty (\hat z)}$. So we have $H^+_{\tilde z}= H^+_{\tilde y}$. This implies that for large $n$, if $\Pi_{\tilde z_n}$ and $\Pi_{\tilde y_n}$ are disjoint then we have either $H^+_{\Pi_{\tilde z_n}}\subsetneq H^+_{\Pi_{\tilde y_n}}$ or $H^+_{\Pi_{\tilde y_n}}\subsetneq H^+_{\Pi_{\tilde z_n}}$. We get a contradiction with the fact that we have $\tilde y_n\subset C(\rho_n)\subset H^+_{\Pi_{\tilde z_n}}$ and $\tilde z_n\subset C(\rho_n)\subset H^+_{\Pi_{\tilde y_n} }$. We deduce that, up extracting a subsequence, $\Pi_{\tilde z_n }$ intersects $\Pi_{\tilde y_n }$. It follows that, for $n$ large enough,  $\hat g_n(\Hp^2)$ intersects $\rho_n (a)\circ\hat g_n(\Hp^2)$. Therefore $\hat g_n(\Hp^2)$ and $\rho_n (a)\circ\hat g_n(\Hp^2)$ coincide. Since  $H^+_{\Pi_{\tilde y}}$ is equal to $ H^+_{\Pi_{\tilde z}}$, the dihedral angle $\theta(\Pi_{\hat z_n},\Pi_{\hat y_n})$ tends to $0$. This implies that the distance between $\tilde z_n$ and $\tilde y_n$ measured on $\hat g_n(\Hp^2)$ tends to $0$. Since  $\hat g_\infty$ is a covering map, there is a neighbourhood ${\cal V}(\hat z)\subset\Hp^2$ of $\hat z$ such that, for $n\in{\overline \N}$ large enough, the map $\hat g_{n|{\cal V}(\hat z)}:{\cal V}(\hat z)\rightarrow \hat g_n({\cal V}(\hat z))$ is a homeomorphism. For $n$ large enough, $\hat y_n$ lies in $\hat g_n({\cal V}(\hat z))$. Hence there is $\hat y'_n \in{\cal V}(\hat z)$ such that $\hat g_n(\hat y'_n)=\tilde y_n$. Since $\hat g_n(\hat y'_n)=\rho_n(a)\circ\hat g_n(\hat y)$ and since $r_n(\Gamma_n)$ is the stabiliser of $\hat g_n(\Hp^2)$ in $\rho_n(\pi_1(M))$, there is $a_n\in\Gamma_n$ such that $a_n\hat y=\hat y'_n$. The point $\hat y'_n$ lies in ${\cal V}(\hat z)$, so $a_n$ moves $\hat y$ a bounded distance. Therefore there is a subsequence such that $(a_n)$ converge to an isometry $a_\infty\in\Gamma_\infty$. Moreover $\hat g_n(\hat y'_n)$ tends to $\hat g_\infty(\hat z)$, hence $(\hat y'_n)$ tend to $\hat z$. Thus we get $a_\infty\hat y=\hat z$. This yields a contradiction with the assumption that $\hat z\not\in\Gamma_\infty\hat y$ and concludes the proof of Lemma \ref{momeo}.
\end{proof}

\indent
Let $F$ be a compact subset of $S$ such that the connected components of $S-F$ are infinite annuli. Since $\Gamma_n$ converges geometrically to $\Gamma_\infty$, there are $\eps>0$ and maps $l_n:F\rightarrow \Hp^2/\Gamma_n$ with the following properties:
\begin{description}
\item - $l_n$ is a homeomorphism onto its image;
\item - $l_n(F)$ is a connected component of the $\eps$-thick part of $\Hp^2/\Gamma_n$;
\item - the induced representations $l_{n*}:\pi_1(F)\rightarrow \Gamma_n$ converge to a faithful representation $l_{\infty*}:\pi_1(F)\rightarrow \Gamma_\infty$;
\item - $g_\infty\circ l_\infty$ coincides with $h_\infty$ on $F\subset\partial M$.
\end{description}

\indent
Let us show that, for large $n$, $g_n\circ l_n:F\rightarrow \partial N(\rho_n)$ is isotopic to $h_n$. Let $\hat F$ be the universal cover of $F$ and let $\hat l_n:\hat F\rightarrow\Hp^2$ be a lift of $l_n$. Since $(l_{n*})$ converges algebraically, we can choose the $\hat l_n$ such that they converge to $\hat l_\infty$ on compact sets. Since $(\hat g_n)$ converges to $\hat g_\infty$, the sequence $(\hat g_n\circ\hat l_n)$ converges to $\hat g_\infty\circ\hat l_\infty$ on compact sets.\\
\indent
Let $x\in int(M)$. Since $\sigma_n$ converges to $\sigma_\infty$, there are diffeomorphisms $u_n:M\rightarrow M$ isotopic to the identity which satisfy the following : for any $k>1$ and $r>0$, there exists $n(k,r)$ such that for $n\geq n(k,r)$, the restriction of $u_n$ to $B(x,r)\subset (M,\sigma_\infty)$ is a $k$-quasi-isometry into its image in $(M,\sigma_\infty)$.
For $n\in\overline{\N}$, set $M_n=\Hp^3/\rho_n(\pi_1(M))$ and let $x_n\in M_n$ be the projection of the origin $o\in\Hp^3$. For $n\in\overline{\N}$, the metric $\sigma_n$ yields an identification between $(int(M),\sigma_n)$ and $M_n$ such that $x$ is identified with $x_n$. Thus we can consider the restriction of $u_n$ to $int(M)$ as an homeomorphism $u_n:M_\infty\rightarrow M_n$.\\
\indent
Let $\tilde u_n:\Hp^3\rightarrow\Hp^3$ be a lift of $u_n$ such that $\tilde u_n(o)=o$. Since $(\rho_n(\pi_1(M)))$ converges geometrically to $\rho_\infty(\pi_1(M))$, $\tilde u_n:\Hp^3\rightarrow\Hp^3$ converges to the identity on compact sets (see \cite{bep}).\\
\indent
Since $(\hat g_n\circ\hat l_n)$ converges to $\hat g_\infty\circ\hat l_\infty$ on compact sets, the sequence $\tilde u_n^{-1}\circ\hat g_n\circ\hat l_n$ converges to $\hat g_\infty\circ\hat l_\infty$ on compact sets. Furthermore, $F$ has a compact fundamental domain. Hence\linebreak $u_n^{-1}\circ g_n\circ l_n:F\rightarrow M_\infty$ converges uniformly to $g_\infty\circ l_\infty$ which coincides with $h_\infty$. It follows that for large $n$, $u_n^{-1}\circ g_n\circ l_n$ is isotopic to $h_{\infty|F}$. Therefore $g_n\circ l_n$ is isotopic to $u_n\circ h_{\infty|F}$. Since $u_n$ is isotopic to the identity, $u_n\circ h_{\infty|F}$ is isotopic to $h_n$. Thus we can change $l_n$ by an isotopy such that $g_n\circ l_n$ coincides with $h_n$ on $F$.\\
\indent
Let $\hat\gamma_n$ be the bending measured geodesic lamination of $\hat g_n$ and let $\gamma_n\in {\cal ML}(\Hp^2/\Gamma_n)$ be its projection. Let $k\subset F$ be an arc such that $l_\infty(k)$ is transverse to $|\gamma_\infty|$ and let $\hat k\subset\hat F$ be a lift of $k$. The arcs $\hat l_n(\hat k)$ converge to $\hat l_\infty(\hat k)$. Since the convex pleated surfaces $(\hat g_n,\Gamma_n,r_n)$ converge to $(\hat g_\infty,\Gamma_\infty,r_\infty)$,  Lemma \ref{plissee} yields $\int_{\hat l_n(\hat k)} d\hat\gamma_n\longrightarrow\int_{\hat l_\infty(\hat k)} d\hat\gamma_\infty$. Since, for any $n\in{\overline \N}$, $g_n\circ l_n$ coincides with $h_n$ on $F$, we have $\int_k d\lambda_n=\int_{\hat l_n(\hat k)} d\hat\gamma_n$ where $\lambda_n$ is the bending measured geodesic lamination of $\sigma_n$. Thus we get $\int_k d\lambda_n\longrightarrow \int_k d\lambda_\infty$.\\
\indent
Doing the same for each component of $\partial M-\lambda_\infty^{(p)}$, we conclude that for any arc\linebreak$k\subset\partial M-\lambda_\infty^{(p)}$ transverse to $|\lambda_\infty|$, we have $\int_k d\lambda_n\longrightarrow \int_k d\lambda_\infty$. It follows that any subsequence of $(\dot\lambda_n)$ contains a subsequence converging in ${\cal ML}(\partial M)/{\cal R}$ and that the limit $\dot\lambda$ differs from $\dot\lambda_\infty$ only on $\lambda^{(p)}_\infty$, namely if $\lambda$ is a representative of $\dot\lambda$ then removing from $\lambda$ the closed leaves which lie in $\lambda^{(p)}_\infty$ yields the same measured geodesic lamination as the one obtained by removing $\lambda^{(p)}_\infty$ from $\lambda_\infty$. Thus, if we show that any leaf of $\lambda^{(p)}_\infty$ is a leaf of $\lambda$ with a weight at least equal to $\pi$, we can conclude that $\dot\lambda_n$ converges to $\dot\lambda_\infty$.

\begin{claim}   \label{feupi}
Let $\lambda$ be a representative of $\dot\lambda$; any leaf $c$ of $\lambda^{(p)}$ is a leaf of $\lambda$ and has a weight greater than or equal to $\pi$ (as a leaf of $\lambda$).
\end{claim}

\begin{proof} 
Let $c$ be a leaf of $|\lambda^{(p)}_\infty|\subset|\lambda_\infty|$. Since $\dot\lambda_n$ converges to $\dot\lambda$ which satisfies $|\dot\lambda|\subset |\lambda_\infty|$, then either $c$ is a leaf of  $\lambda$ and we will denote its weight by $w(c )$, or $c$ is a simple closed curve disjoint from $\lambda$ and we will take $w(c)=0$.\\
\indent
We will prove the claim by assuming that $w(c)<\pi$ and by ending in a contradiction.\\
\indent
Let $S$ be a component of $\partial M-\lambda^{(p)}_\infty$ whose closure contains $c$. From the fact that $r_n(\pi_1(S))$ converges geometrically to $r_\infty(\pi_1(S))$ we deduce that $l_{s_n}(c )\longrightarrow 0$.\\
\indent
In the case where, up to extracting a subsequence, $i(c,\lambda_n)$ is equal to $0$ for all $n$, it follows from \cite[Lemma 19]{meister1} that $c$ is a leaf of $\lambda$ with a weight equal to $\pi$. Now we can assume that we have $i(c,\lambda_n)>0$ for any $n$ large enough. Especially $c$ is not a leaf of $\lambda_n$.\\
\indent
Let $\eta>0$ be a number such that $\pi-\eta$ is greater than $w(c  )$. The curve $c  $ lies in the boundaries of two surfaces $F^1$ and $F^2$ (which may coincide) such that the following holds : $c\subset int(F^1\cup F^2)$ and either $F^1$ (resp. $F^2$) is a pair of pants satisfying $int(F^1)\cap |\lambda|=\emptyset$ or there is connected component of $\lambda$ which is an arational geodesic lamination in $F^1$ (resp. $F^2$). We will only deal with the case where $F^1$ and $F^2$ are distinct, the other case is handled in the same way. Approximating $|\lambda|\cap F_i$ by an arc with endpoints in $c$, we can construct a simple closed geodesic $d$ that intersects $c$ in two points and satisfies $i(d,\lambda-c  )\leq\frac{\eta}{8}$. Consider a point $x_n^1$ (respectively $x_n^2$) of $d\cap F^1-\lambda_n$ (respectively  $d\cap F^2-\lambda_n$) lying in the thick part of $(F^1,s_n)$ (respectively $(F^2,s_n)$). Let $k_n^1$ and $k_n^2$ be the connected components of $d-\{x_n^1,x_n^2\}$. We have $\limsup \int_{k_n^i} d\lambda_n\leq i(d,\lambda-c)+w(c)$.\\
\indent
Now we want to find $2$ points $y_n^1\subset k_n^1$ and $y_n^2\subset k_n^2$ that cut $k_n^1$ and $k_n^2$ into $4$ arcs $\kappa^j_n$,  $1\leq j\leq 4$, such that we have $\limsup\int_{\kappa^j_n}\! d\lambda_n<\frac{\pi}{2}$ for any $j$. If $w(c)=0$ then we have $\limsup\int_{k_n^i} d\lambda_n\leq i(d,\lambda)=i(d,\lambda-c  )\leq\frac{\eta}{4}<\frac{\pi}{4}$. Hence any points $y_n^1\subset k_n^1$ and $y_n^2\subset k_n^2$ would be suitable. If we have $w(c)>0$ then $c$ is a leaf of $\lambda$. Since we have assumed that $c$ is not a leaf of $\lambda_n$, $\lambda_n$ spirals more and more around $c$. More precisely, in a neighbourhood ${\cal V}(c)$ of $c$, $\lambda_n\cap {\cal V}(c)$ is a family of compact arcs each one spiralling many times toward $c$ and carrying a small measure. It follows that we can find points $y_n^j\subset k_n^j$ such that if $\kappa_n^j$, $1\leq j\leq 4$, are the closure of the components of $d-\{x_n^1,x_n^2,y_n^1,y_n^2\}$ (cf. figure \ref{fig3}), we have $\frac{w(c)}{2}-\frac{\eta}{8}\leq\limsup\int_{\kappa_n^j} d\lambda_n\leq i(d,\lambda-c  )+\frac{w(c)}{2}+\frac{eta}{8}\leq\frac{\pi}{2}-\frac{\eta}{4}$. Roughly, the $y_n^j$ are chosen so that $\lambda_n$ spirals as many time before $y_n^j$ as it does after.
\begin{figure}[hbtp]
\psfrag{a}{$\kappa_n^1$}
\psfrag{b}{$\kappa_n^2$}
\psfrag{c}{$\kappa_n^3$}
\psfrag{d}{$\kappa_n^4$}
\psfrag{e}{$x_n^1$}
\psfrag{f}{$x_n^2$}
\psfrag{g}{$y_n^1$}
\psfrag{h}{$y_n^2$}
\psfrag{i}{$\lambda_n$}
\centerline{\includegraphics{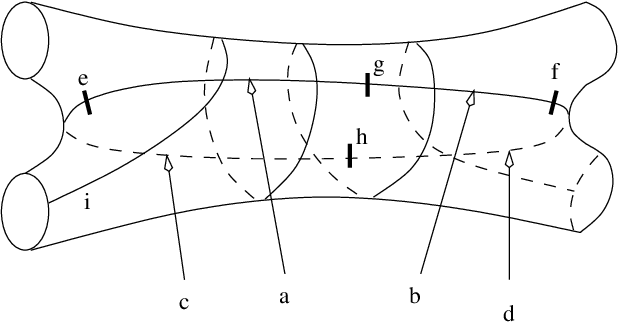}}
\caption{Cutting $d$}
\label{fig3}
\end{figure}
For $j\in\{1,2,3,4\}$, consider a lift $\tilde\kappa_n^j\subset C(\rho_n)\subset\Hp^3$ of $\kappa_n^j$ such that $\bigcup_j\tilde\kappa_n^j$ is connected. Let $\tilde d_n^j\subset\Hp^3$ be the geodesic segment joining the vertices of $\tilde\kappa_n^j$.\\
\indent
If there is $j_0$ such that $l_{s_n}(\kappa_n^{j_0})$ is bounded, then $\kappa_n^{j_0}$ is entirely contained in the thick part of $(F^1\cup F^2,s_n)$. Since $g_n(F^1)$ and $f_n(F^2)$ tend to pleated surfaces whose bending measured geodesic laminations are $\lambda\cap F^j $, we have the following inequality $\lim\int_{\kappa_n^{j_0}} d\lambda_n\leq i(d,\lambda-c  )\leq\frac{\eta}{8}$. Hence we have $\frac{w(c)}{2}\leq\frac{\eta}{4}$. Thus we get $i(d,\lambda)\leq 2w(c)+\frac{\eta}{4}\leq 2\eta$. Since $d$ intersects $c$ transversally, and since $l_{s_n}(c)\longrightarrow 0$, $l_{s_n}(d)$ tends to $\infty$. Taking $\eta< \frac{\pi}{4}$, we get from the slight bending lemmas (Lemma \ref{approx}) that $l_{\sigma_n}(d^*)\geq C_{2\eta}(l_{s_n}(d)+A_{5\eta})$. Thus $l_{\sigma_n}(d^*)$ tends to $\infty$ contradicting the algebraic convergence of $\rho_n$.\\
\indent
Thus we have $\forall j$, $l_{s_n}(\kappa_n^j)\longrightarrow\infty$. The slight bending Lemmas (Lemma \ref{lemap}) says that there is $C>0$ such that $l_{s_n}(\tilde\kappa_n^j)\leq C l_{\sigma_n}(\tilde d_n^j)$ and that the interior angle between two adjacent segments $\tilde d_n^j$ is greater than $\frac{\eta}{4}$. Let $a$ be the element of $\pi_1(M)$ such that $\rho_n(a)$ fixes $d^*$ and let $\tilde d_n=\bigcup_{i\in\Z;j=1,2,3,4} a^i(\kappa_n^j)$. The curve $\tilde d_n$ is the union of long geodesic segments with interior angles greater than $\frac{\eta}{4}$. It follows from a classical result (see \cite{meister2} for example) that there is $K$ such that $l(\tilde d_n/a)$ is smaller than $K l_{\sigma_n}(d^*)$. Thus we get $l_{s_n}(d)\leq C l(\tilde d_n/a)\leq CK l_{\sigma_n}(d^*)$. Since $l_{s_n}(d)\longrightarrow \infty$, this yields a contradiction with the algebraic convergence of $\sigma_n$.
\end{proof}

\indent
 It follows that $\dot\lambda=\dot\lambda_{\infty}$. Thus Proposition \ref{con} is proved.
\end{proof}

\newpage

\section{Erratum, octobre 2025}
We need to correct an error in the present article. As stated the main theorem (Theorem 1) is false. In the introduction, we equip the space $\ML(\partial M)/\CR$ with the quotient topology but the bending map is not continuous with respect to this topology (see comment at the end of this note). On the other hand all the arguments in \cite{lecuire:continuity} are correct if we use another topology on $\ML(\partial M)/\CR$ that we will call the {\em tubular topology} and define below. \\


A first step in defining the tubular topology consists in constructing a specific basis of neighborhoods for the weak$^*$ topology on $\ML(S)$ based on a construction due to Otal, \cite[\textsection 3.2]{otal:fibre}. Let $\lambda$ be a measured geodesic lamination on $S$ and let $L$ be a connected component of its support $|\lambda|$. If $L$ is a simple closed curve, we choose a geodesic arc that intersects $L$ once and is disjoint from the other components of $L$. If $L$ is not a simple closed curve, we start with a geodesic arc $k$ that intersect $L$ such that $\int_{k} d\lambda<\frac{\pi}{2}$. By \cite[Proposition A.3.4]{otal:fibre}, every leaf of $L$ is dense in $L$ and by \cite[\textsection 3.2]{otal:fibre} the closures of the components of $L\setminus k$ form $p_L$ families such that any two arcs in a family are isotopic relative to $k$. As explained in \cite[\textsection 3.2]{otal:fibre}, this allows us to build $p_L$ rectangles $r_i:[0,1]\times[0,1]\to S$ such that $r_i(\{\eta\}\times[0,1])\subset k$ for $\eta\in\{0,1\}$, that $r_i([0,1]\times\{\eta\})$ are geodesic arcs for $\eta\in\{0,1\}$, that the restriction of $r_i$ to $(0,1)\times[0,1]$ is an embedding, that $r_i((0,1)\times[0,1])\cap r_j((0,1)\times[0,1])=\emptyset$ for $i\neq j$ and that $L\subset\bigcup_i r_i([0,1]\times[0,1])$. We get a family of arcs by taking the sides $r_i(\{\eta\}\times[0,1])\subset k$ and $r_i([0,1]\times\{\eta\})$, $\eta\in\{0,1\}$, for all $r_i$. We do this construction for each component of $|\lambda|$ and obtain a family of arcs $k_i, i\leq p$ such that $\int_{k_i} d\lambda<\frac{\pi}{2}$ unless $k_i$ intersects a closed leaf of $\lambda$. We add a geodesic arc in each component of $S\setminus|\lambda|$ to obtain a family of arcs $k_i, i\leq q$ such that the intersection between $\bigcup_i k_i$ and any simple geodesic is transverse and non-empty (we may need to slightly change the arcs $k_i$, compare with \cite[p.19]{bonahon:laminations}).
Set $\CB_\eps(\lambda)=\{\mu\in\ML(S)\mbox{ such that } |\int_{k_i}d\mu-\int_{k_i} d\lambda|<\eps \mbox{ for any } i\leq q\}$. It follows from the work of Thurston \cite[\textsection 8.2]{thurston:notes} that $\{\CB_\eps(\lambda)|\eps>0\}$ is a local basis at $\lambda$.\\

Given $\dot\lambda\in\ML(S)/\CR$, we denote by $\lambda'\in\ML(S)$ the representative of $\dot\lambda$ that has no leaf with a weight greater than $\pi$ and by $F_{\dot\lambda}$ be the subset of $\ML(S)$ made up of all measured geodesic laminations that project to $\dot\lambda$. Let $\lambda^{(p)}$ be the union of the leaves of $\lambda'$ with a weight equal to $\pi$. Using the previous construction, we define the set $\mathring\CB_\eps(\dot\lambda,\eta)=\{\mu\in\ML(S)$ such that $\int_{k_i}d\mu>\pi-\eps$ if $k_i$ intersects $\lambda^{(p)}$ and $|\int_{k_i}d\mu-\int_{k_i} d\lambda'|<\eps$ otherwise$\}$ which is an open set containing $F_{\dot\lambda}$.  Let $m(\dot\lambda)$ be the maximum of the weights of the closed leaves of $\lambda'$ with weights smaller than $\pi$. For $\eps<\min\{\frac{\pi}{2},\pi-m(\lambda)\}$, any leaf of a measured geodesic lamination $\mu\in \mathring\CB_\eps(\dot\lambda)$ with a weight larger than or equal to $\pi$ is a leaf of $\lambda^{(p)}$. It follows that $\mathring\CB_\eps(\dot\lambda)$ is saturated with respect to $\CR$. Let $\dot\CB_\eps(\dot\lambda)$ be the projection of $\mathring\CB_\eps(\dot\lambda)$ to $\ML(S)/\CR$. We define the tubular topology by setting that $\{\dot\CB_\eps(\dot\lambda)|\eps>0\}$ is a local basis at $\dot\lambda$.\\

Since the set $\{\dot\CB_{\frac{1}{n}}(\dot\lambda)|n\in\N\}$ is a countable local basis, $\ML(S)/\CR$ equipped with the tubular topology is first countable and its topology is determined by its converging sequences. Thus we can analyze the continuity of the bending map simply by studying the behavior of images of converging sequences. In \cite{lecuire:properness}, it is also proved that, with respect to the tubular topology, compactness and sequential compactness are equivalent.\\

We can now correct the main statement :

\begin{theoremc}       \label{three}
	The map $b_{\CR}$ from ${\mathcal GF}(M)$ to ${\ML}(\partial M)/{\CR}$ equipped with the tubular topology is a continuous map.
\end{theoremc}

The arguments developed in present paper actually give a complete proof of this statement since they only use the following characterization of the topology on ${\ML}(S)/{\CR}$: a sequence $(\dot\lambda_n)\subset {\ML}(S)/{\CR}$ converges to $\dot\lambda\in {\ML}(S)/{\CR}$ if and only if $\int_k d\lambda'_n$ converges to $\int_k d\lambda'$ for any arc $k$ that is disjoint from $\lambda^{(p)}$ and $\liminf\int_k d\lambda'_n\geq\pi$ for any arc $k$ that intersects $\lambda^{(p)}$. This property characterizes the tubular topology, this follows from the definition and the work of Thurston (\cite[\textsection 8 and 9]{thurston:notes}, see also \cite{penner:harer} and \cite{otal:fibre}), but not the quotient topology as illustrated below.\\

To conclude this erratum, we will show that the quotient topology on $\ML(S)/\CR$ is different from the tubular topology and that the bending map is not continuous when ${\ML}(\partial M)/{\CR}$ is equipped with the quotient topology. 

Let $c,d\subset S$ be two simple closed curves that intersect once and let $D_c:S\to S$ be the right Dehn twist about $c$. Let $\gamma_n\in\ML(S)$ be obtained by equipping the geodesic $d_n$ in the isotopy class of $D_c^{n^3}(d)$ with a transverse Dirac measure with weight $\frac{1}{n}$ and let $\gamma\in\ML(S)$ be $c$ with weight $\pi$. Let $k_d\subset S\setminus c$ be an arc intersecting each $d_n$ once and $k_c\subset S\setminus d$ be an arc intersecting $c$ once. Since $\int_{k_d} d\gamma_n=\frac{1}{n}$ and $\int_{k_c} d\gamma_n\sim\frac{n^3}{n}$, it is not difficult to prove that, for any $\eps>0$, $\mathring{\CB}_\eps(\dot\gamma)$ contains $\gamma_n$ for $n$ large enough (depending on $\eps$)  and that $(\dot\gamma_n)$ converges to $\dot\gamma$ for the tubular topology. On the other hand, the set $\mathcal V=\{\mu\in\ML(S)$ such that $\int_{k_d}d\mu<(\int_{k_c}d\mu)^{-1}\}$ is an open subset of $\ML(S)$ that contains $\gamma$ but not $\gamma_n$ for $n$ large. The intersection $\mathcal U=\mathcal V\cap \mathring{\CB}_\eps(\dot\gamma)$ is an open set which is saturated with respect to $\CR$ when $\eps<\{\frac{\pi}{2}\}$. Hence it projects to an open subset $\dot{\mathcal U}$ of $\ML(S)/\CR$ for the quotient topology. Since $\dot{\mathcal U}$ contains $\dot\gamma$ but not $\dot\gamma_n$ for $n$ large, $(\dot\gamma_n)$ does not converge in the quotient topology.\\

Let $M$ be a compact acylindrical hyperbolic $3$-manifold whose boundary is homeomorphic to $S$. For any $n$, there is a convex cocompact metric $\sigma_n$ with bending measured geodesic lamination $\gamma_n$. By \cite{thurston:hypI}, $\sigma_n$ has a converging subsequence. The fact that $(\dot\gamma_n)$ converges in the tubular topology but not in the quotient topology shows that the corrected version of Theorem 1 contradict the original version. 

\end{document}